\documentclass[11pt,leqno]{article}
\usepackage{amsfonts}
\usepackage{amsmath}
\usepackage{amssymb}
\newcommand{\pref}[1]{(\ref{#1})}
\newtheorem{theo}{Theorem}[section]
\newtheorem{lem}[theo]{Lemma}
\newtheorem{prop}[theo]{Proposition}
\newtheorem{cor}[theo]{Corollary}

\title{\Large \bf {$Sp(n)U(1)$-connections with parallel totally skew-symmetric torsion}}
\author{{\sc Bogdan Alexandrov}
\thanks{Supported by SFB 288 "Differential geometry and quantum
physics" of DFG	and The European Contract Human Potential Programme,
Research Training Network HPRN-CT-2000-00101}
}
\date{}
\begin{document}
\maketitle
\vspace{5mm}
\begin{abstract}
We consider the unique Hermitian connection with totally skew-symmetric torsion on a Hermitian manifold. We prove that if the torsion is parallel and the holonomy is $Sp(n)U(1) \subset U(2n) \times U(1)$, then the manifold is locally isomorphic to the twistor space of a quaternionic K\"ahler manifold with positive scalar curvature. If the manifold is complete, then it is globally isomorphic to such a twistor space.
\\[10mm]
\\[10mm]
\end{abstract}

\section{Introduction}

Let $(M^n,g)$ be an oriented Riemannian manifold. Let $G$ be a subgroup of $SO(n)$ and $P_G$ be a $G$-structure on $M$, i.e., $P_G$ is a principal $G$-bundle which is subbundle of the bundle of oriented orthonormal frames $P_{SO(n)}$. Suppose that the Levi-Civita connection $\nabla$ does not come from a connection on $P_G$. Which is the best connection on $P_G$ in this case?

The first obvious choice is the canonical connection. It is the unique connection $\nabla^c$ whose torsion is the intrinsic torsion of the $G$-structure $P_G$. It can be thought as the orthogonal projection of the Levi-Civita connection in the affine space of all $G$-connections in the following sense: $\nabla^c = \nabla + A^c$, where at each point $p \in M$ $A^c_p$ is orthogonal to $(T_p M)^* \otimes \mathfrak{g}$ in $(T_p M)^* \otimes so(T_p M)$.

Another choice would be to replace the condition of vanishing torsion, which characterizes the Levi-Civita connection, by the requirement that the torsion is (covariantly) constant. This implies the existence of an invariant element of $(T_p M)^* \otimes so(T_p M)$ with respect to the holonomy group of the connection. If there is no such invariant element with respect to $G$ itself, this would mean that the holonomy group is a proper subgroup of $G$, i.e., a further reduction of the structure group should be possible. Thus a $G$-connection with parallel torsion does not always exist.

A third possibility is a $G$-connection $\nabla^a = \nabla + A^a$, for which the potential $A^a$ (or, equivalently, the torsion $T^a$) is totally skew-symmetric. The advantage of such a connection is that it has the same geodesics as the Levi-Civita connection. In particular, it is complete exactly when the metric is complete. In the general case though there is neither existence nor uniqueness of a $G$-connection with totally skew-symmetric torsion.

As an illustration let us consider the case of almost Hermitian structure, i.e., $G = U(n) \subset SO(2n)$. It is proved in \cite{FI} (see also \cite{G}) that a $U(n)$-connection (or, in the established terminology, a Hermitian connection) $\nabla^a$ with totally skew-symmetric torsion exists iff the Nijenhuis tensor is totally skew-symmetric and in this case it is unique. The last condition means that the almost Hermitian manifold lies in the Gray-Hervella class $\mathcal{G}_1$, which contains the nearly K\"ahler and the Hermitian manifolds \cite{GH}. The connection $\nabla^a$ has been used by Bismut to prove a local index theorem for the Dolbeault operator on Hermitian non-K\"ahler manifolds \cite{Bis}.  The perfect situation occurs in the complementary case of nearly K\"ahler manifolds: the canonical connection has totally skew-symmetric torsion, which is furthermore parallel. The last result was proved by Kirichenko \cite{K} (for another proof see \cite{BM}). On the other hand, on a manifold of class $\mathcal{G}_1$ which is not nearly K\"ahler (in particular, on a Hermitian manifold) neither $\nabla^a$ coincides with the canonical connection nor is its torsion $T^a$ parallel in general.

In this paper we are interested in $Sp(n)U(1)$-structures, where $Sp(n)U(1)$ is considered as a subgroup of $U(2n) \times U(1) \subset U(2n+1) \subset SO(4n+2)$ by certain inclusion $\rho$ (see section~\ref{sec1}). More precisely, we study the $(4n+2)$-dimensional Hermitian manifolds whose unique Hermitian connection with totally skew-symmetric torsion $\nabla^a$  has holonomy contained in $\rho (Sp(n)U(1))$ and parallel torsion $T^a$.

In section~\ref{sec2} we consider $(2m+2)$-dimensional Hermitian manifolds such that the holonomy group of $\nabla^a$ is a subgroup of $U(m) \times U(1)$ and the torsion  $T^a$ satisfies a simple algebraic condition (see Proposition~\ref{prop4}). We show that there is an interesting correspondence between such manifolds and nearly K\"ahler manifolds with $Hol(\nabla^c) \subset U(m) \times U(1)$. This allows us to prove that the torsion $T^a$ is parallel. Furthermore, we prove that if $T^a$ is non-degenerate, then $m=2n$ and $Hol(\nabla^a) \subset \rho (Sp(n)U(1))$. It follows also that the Ricci tensors of both $\nabla^a$ and the Levi-Civita connection $\nabla$ are $\nabla^a$-parallel and positive definite.

In section~\ref{sec5} we consider the curvature tensor of a Hermitian manifold such that $\nabla^a$ has parallel torsion and holonomy $\rho (Sp(n)U(1))$. It turns out that it decomposes in a way very similar to the decomposition of the curvature of a quaternionic K\"ahler manifold. The methods of this section can be applied to give a proof of the above mentioned fact that the torsion of a nearly K\"ahler manifold is parallel. This is done using the first Bianchi identity and the representation theory of $U(n)$.

In section~\ref{sec3} we show what section~\ref{sec5} makes conceivable: there are examples on twistor spaces of quaternionic K\"ahler manifolds (we adopt the definition that a 4-dimensional manifold is quaternionic K\"ahler if it is self-dual and Einstein). The twistor space $\mathcal{Z}$ of a quaternionic K\"ahler manifold $M'$ carries in a natural way two almost complex structures $J_1$ (integrable) and $J_2$ (non-integrable). They have been first defined on twistor spaces of 4-dimensional manifolds in \cite{AHS} and \cite{ES} respectively. On $\mathcal{Z}$ exists also a one-parameter family of metrics $h_t$, $t>0$, Hermitian with respect to both $J_1$ and $J_2$. If the base $M'$ has positive scalar curvature, there are two particularly interesting values $t_0$ and $t_1$ of the parameter $t$, such that $(\mathcal{Z},h_{t_0},J_1)$ is K\"ahler and $(\mathcal{Z},h_{t_1},J_2)$ is nearly K\"ahler \cite{FK,S,M,AGI}. We prove that the connection $\nabla^a$ of $(\mathcal{Z},h_{t_1},J_1)$ has parallel torsion and holonomy $\rho (Sp(n)U(1))$.

In the next two sections we prove the main result of this paper: there are no other examples, both globally and locally.

The global result is contained in Theorem~\ref{thm1}: a complete Hermitian manifold, such that $\nabla^a$ has parallel torsion and holonomy $\rho (Sp(n)U(1))$, is the twistor space $(\mathcal{Z},h_{t_1},J_1)$ of some compact quaternionic K\"ahler manifold $M'$ with positive scalar curvature. In particular, there are only finitely many such manifolds in each dimension since the same is true for the compact quaternionic K\"ahler manifolds with positive scalar curvature \cite{L,LS}. In fact, the only known examples of compact quaternionic K\"ahler manifolds with positive scalar curvature are the Wolf spaces \cite{W,B}. In this case $(\mathcal{Z},h_{t_1},J_1)$ is homogeneous naturally reductive. According to the results in \cite{FK,H,PS,HH} the Wolf spaces are the only compact quaternionic K\"ahler manifolds with positive scalar curvature in dimensions 4, 8 and 12. Therefore we get a complete list of the Hermitian manifolds of dimension 6, 10 and 14, satisfying the above conditions. In the proof of Theorem~\ref{thm1} we use the corresponding results about nearly K\"ahler manifolds of Belgun--Moroianu \cite{BM} and Nagy \cite{N}.

Section~\ref{sec6} is devoted to the proof of Theorem~\ref{thm2}, which is the local version of Theorem~\ref{thm1}. As a corollary we get a corresponding local result for nearly K\"ahler manifolds.

In the Appendix we gather some simple and definitely well known facts about projections of tensors, subbundles of tensor bundles and connections by submersions. Some of them could be found for example in \cite{IK} but the author was unable to find others in explicit form in the literature (especially those about projections of subbundles).

Finally, we should mention that in dimension 6 this paper covers one of the cases considered in \cite{AFS}. The subject of \cite{AFS} are the 6-dimensional Hermitian manifolds on which $\nabla^a$ has parallel torsion $T^a$. Thus $T^a$ is invariant with respect to the holonomy group of $\nabla^a$ and this strongly restricts the possible holonomy groups. In \cite{AFS} these possibilities are listed and the corresponding manifolds are studied.

\section{Algebraic preliminaries}\label{sec1}

Let $T \cong \mathbb{R}^{2m+2}$ be the standard $(2m+2)$-dimensional real representation of $U(m+1)$ (in the following sections $T$ will be the tangent space of a Hermitian manifold). Its complexification is $T^{\mathbb{C}} = T^{1,0} \oplus T^{0,1}$, where $T^{1,0} \cong \mathbb{C}^{m+1}$ is the standard $(m+1)$-dimensional complex representation of $U(m+1)$ and $T^{0,1}$ is its conjugate. Denote by $\Lambda ^{p,q} T^* \cong \Lambda ^p (T^{1,0})^* \otimes \Lambda ^q (T^{0,1})^*$ the space of (complex) $(p,q)$-forms on $T$. Let $\mathcal{H} \cong \mathbb{R}^{2m}$ (resp. $\mathcal{V} \cong \mathbb{R}^2$) be the standard $2m$-dimensional (resp. 2-dimensional) real representation of $U(m)$ (resp. $U(1)$), with further notations similar as above for $T$. Then, as representations of $U(m) \times U(1)$,
$$T = \mathcal{H} \oplus \mathcal{V}, \quad T^{1,0} = \mathcal{H}^{1,0} \oplus \mathcal{V}^{1,0}, \quad T^{0,1} = \mathcal{H} ^{0,1}\oplus \mathcal{V}^{0,1}.$$

Let us now consider the subgroup $Sp(n)U(1)$ of $U(2n)$ (as an abstract group it is isomorphic to $(Sp(n) \times U(1))/_{\mathbb{Z}_2}$). Define the inclusion $\rho : Sp(n)U(1) \longrightarrow U(2n) \times U(1) \subset U(2n+1)$ as follows:
if $a \in Sp(n)$, $b \in U(1)$, then
$$\rho (ab) =
\begin{pmatrix}
ab & 0 \\
0 & b^2
\end{pmatrix}
\in U(2n) \times U(1) \subset U(2n+1).$$
If $n=1$, $Sp(1)U(1) \cong U(2)$ and for $c \in U(2)$
$$\rho (c) =
\begin{pmatrix}
c & 0 \\
0 & \det c
\end{pmatrix}
\in U(2) \times U(1) \subset U(3).$$

Denote by $E$ the standard $2n$-dimensional complex representation of $Sp(n)$ and by $F(k)$ the complex (1-dimensional) representation of $U(1)$ with weight $k$. Since $E$ is self-adjoint, by the definition of $\rho$ we have that, as $\rho (Sp(n)U(1))$-representations,
$$\mathcal{H}^{1,0} \cong (\mathcal{H} ^{0,1})^* \cong E \otimes F(1), \quad \mathcal{H}^{0,1} \cong (\mathcal{H} ^{1,0})^* \cong E \otimes F(-1),$$
$$\mathcal{V}^{1,0} \cong (\mathcal{V} ^{0,1})^* \cong F(2), \quad \mathcal{V}^{0,1} \cong (\mathcal{V} ^{1,0})^* \cong F(-2).$$
Let $e_1,\dots ,e_{2n+1}$ be an orthonormal basis of $T^{1,0}$ such that $e_1,\dots ,e_{2n} \in \mathcal{H}^{1,0}$, $e_{2n+1} \in \mathcal{V}^{1,0}$ and
$$\omega_0 = \sum _{k=1}^n e^{2k-1} \wedge e^{2k} \in \Lambda ^{2,0} \mathcal{H}^* \cong \Lambda^2 E \otimes F(-2)$$
is the 2-form corresponding to the $Sp(n)$-invariant 2-form on $\Lambda^2 E$ ($e^1,\dots ,e^{2n+1}$ is the basis dual to $e_1,\dots ,e_{2n+1}$).

Denote
$$T_0 :=\omega_0 \wedge \bar{e}^{2n+1} \in \Lambda ^{2,0} \mathcal{H}^* \otimes \Lambda ^{0,1} \mathcal{V}^* \subset \Lambda ^{2,1} T^*.$$
$T_0$ is obviously $\rho (Sp(n)U(1))$-invariant and
\begin{equation}\label{101}
|\lambda T_0 + \bar{\lambda} \overline{T}_0|^2 = 12n|\lambda|^2.
\end{equation}
(Here and in the sequel we use the tensorial norms. The other popular convention is that the norm of a (skew-symmetric) $k$-form is its tensorial norm divided by $k!$).

\begin{prop}\label{prop1}
The subspace of $\Lambda ^{2,1} T^*$ on which $\rho (Sp(n)U(1))$ acts trivially is 1-dim\-en\-sion\-al and is spanned by $T_0$.
\end{prop}

\noindent {\it Proof:}
We have
$$\Lambda ^{2,1} T^* \cong \Lambda ^{2,0} T^* \otimes (T^{0,1})^* \cong (\Lambda ^{2,0} \mathcal{H}^* \oplus (\mathcal{H}^{1,0})^* \otimes (\mathcal{V}^{1,0})^*) \otimes ((\mathcal{H}^{0,1})^* \oplus (\mathcal{V}^{0,1})^*).$$
Further, the decompositions of these tensor products into irreducible $\rho (Sp(n)U(1))$-re\-pre\-sen\-ta\-tions are
$$\Lambda ^{2,0} \mathcal{H}^* \otimes (\mathcal{H}^{0,1})^* \cong ((\mathbb{C} \oplus \Lambda^2_0 E) \otimes F(-2)) \otimes (E \otimes F(1)) \cong (E \oplus E \oplus \Lambda^3_0 E \oplus K) \otimes F(-1),$$
where $K$ is the irreducible representation of $Sp(n)$ with highest weight $(2,1,0,\dots ,0)$,
$$\Lambda ^{2,0} \mathcal{H}^* \otimes (\mathcal{V}^{0,1})^* \cong ((\mathbb{C} \oplus \Lambda^2_0 E) \otimes F(-2)) \otimes F(2) \cong \mathbb{C} \oplus \Lambda^2_0 E,$$
$$((\mathcal{H}^{1,0})^* \otimes (\mathcal{V}^{1,0})^*) \otimes (\mathcal{H}^{0,1})^* \cong (E \otimes F(-1) \otimes F(-2)) \otimes (E \otimes F(1)) \cong (\mathbb{C} \oplus \Lambda^2_0 E \oplus S^2 E) \otimes F(-2),$$
$$((\mathcal{H}^{1,0})^* \otimes (\mathcal{V}^{1,0})^*) \otimes (\mathcal{V}^{0,1})^* \cong (E \otimes F(-1) \otimes F(-2)) \otimes F(2) \cong E \otimes F(-1)$$
($\Lambda^2_0 E$ and $K$ are 0 if $n=1$ and the same is true for $\Lambda^3_0 E$ if $n \leq 2$).

Hence the only subspace of $\Lambda ^{2,1} T^*$, on which $\rho (Sp(n)U(1))$ acts trivially,  is contained in $\Lambda ^{2,0} \mathcal{H}^* \otimes (\mathcal{V}^{0,1})^*$ and is obviously spanned by $T_0$. \hfill $\Box $

\begin{cor}\label{cor1}
Consider the space of real $((2,1) \oplus (1,2))$-forms as a real $\rho (Sp(n)U(1))$-representation. Then the subspace on which $\rho (Sp(n)U(1))$ acts trivially has (real) dimension 2 and is spanned by $T_0 + \overline{T}_0$ and $i T_0 - i \overline{T}_0$.
\end{cor}

\begin{prop}\label{prop2}
The subgroup of $U(2n+1)$, which preserves $T_0$, is $\rho (Sp(n)U(1))$.
\end{prop}

\noindent {\it Proof:}
We already know that $T_0$ is $\rho (Sp(n)U(1))$-invariant. Let $h \in U(2n+1)$ be such that $h(T_0) = T_0$, i.e.,
\begin{equation}\label{2001}
h(\omega_0) \wedge h(\bar{e}^{2n+1}) = \omega_0 \wedge \bar{e}^{2n+1}.
\end{equation}
Hence $h$ preserves $(\mathcal{V}^{0,1})^* = span \{ \bar{e}^{2n+1} \} $ since the splitting $T^{\mathbb{C}} = T^{1,0} \oplus T^{0,1}$ is $U(2n+1)$-invariant. Therefore
$$ h=
\begin{pmatrix}
c & 0 \\
0 & d
\end{pmatrix}
\in U(2n) \times U(1).$$
Let $a \in U(2n)$, $b \in U(1)$ be such that $c=ab$, $d=b^2$. Then $h(\bar{e}^{2n+1}) = b^2.\bar{e}^{2n+1}$, $h(\omega_0) = b^{-2}.a(\omega_0)$ and so $h(\omega_0) \wedge h(\bar{e}^{2n+1}) = a(\omega_0) \wedge \bar{e}^{2n+1}$. This and \pref{2001} imply $a(\omega_0) = \omega_0$, i.e., $a \in Sp(n)$ and therefore $h \in \rho (Sp(n)U(1))$. \hfill $\Box$

\vspace{3mm}
Let $g$ denote the standard inner product on $T$ and $J$ be the standard complex structure (acting as multiplication by $i$ on $T \cong \mathbb{C}^{2n+1}$). Define $K,I \in End(T)$ by
\begin{equation}\label{102}
g((K+iI)X,Y) = 2 \omega_0 (X,Y), \quad X,Y \in T.
\end{equation}
Then $K$ and $I$ vanish on $\mathcal{V}$, preserve $\mathcal{H}$, $K_{|\mathcal{H}}$ and $I_{|\mathcal{H}}$ are orthogonal with respect to $g_{|\mathcal{H}}$ and $I_{|\mathcal{H}}$, $J_{|\mathcal{H}}$, $K_{|\mathcal{H}}$ satisfy the quaternionic identities. Since $span_{\mathbb{C}} \{ \omega_0 \}$ is $\rho (Sp(n)U(1))$-invariant, \hfill $span \{ K,I \}$ \hfill is \hfill $\rho (Sp(n)U(1))$-invariant. \hfill Thus \hfill $span \{ I_{|\mathcal{H}}, J_{|\mathcal{H}}, K_{|\mathcal{H}} \}$ \hfill is \hfill a \\
$\rho (Sp(n)U(1))$-invariant quaternionic structure on $\mathcal{H}$ compatible with $g_{|\mathcal{H}}$.

Later we shall need also a second inclusion $\rho_2 : Sp(n)U(1) \longrightarrow U(2n) \times U(1) \subset U(2n+1)$, defined as follows:
if $a \in Sp(n)$, $b \in U(1)$, then
$$\rho_2 (ab) =
\begin{pmatrix}
ab & 0 \\
0 & b^{-2}
\end{pmatrix}
\in U(2n) \times U(1)$$
(if $n=1$, we have $\rho_2 : U(2) \longrightarrow U(2) \times U(1)$:
$$\rho_2(c) =
\begin{pmatrix}
c & 0 \\
0 & (\det c)^{-1}
\end{pmatrix}
).$$

\section{Hermitian manifolds with $Hol(\nabla^a) \subset U(m) \times U(1)$}\label{sec2}

Let $(M,g,J)$ be an almost Hermitian manifold. We denote the Levi-Civita connection by $\nabla$ and the K\"ahler form by $\Omega$,
$$\Omega (X,Y) = g(JX,Y).$$
We shall assume that $(M,g,J)$ belongs to the class $\mathcal{G}_1$ of Gray-Hervella \cite{GH}. This class is characterized by the property that the Nijenhuis tensor $N$ is totally skew-symmetric. As shown in \cite{FI}, this is equivalent to the existence of a Hermitian connection $\nabla^a$ with totally skew-symmetric torsion $T^a$ and this connection is furthermore unique. It is given by
$$\nabla^a = \nabla +\frac{1}{2} T^a,$$
$$T^a = -d^c \Omega + N,$$
where $d^c \Omega (X,Y,Z) = -d \Omega (JX,JY,JZ)$. The class $\mathcal{G}_1$ contains as subclasses the Hermitian manifolds and the nearly K\"ahler manifolds. They can be distinguished in terms of the torsion $T^a$ as follows: the manifold is Hermitian (i.e., $J$ is integrable) iff $T^a$ is a $((2,1) \oplus (1,2))$-form and in this case $T^a = -d^c \Omega$, and it is nearly K\"ahler iff $T^a$ is a $((3,0) \oplus (0,3))$-form. In the latter case $\nabla^a$ coincides with the canonical Hermitian connection $\nabla^c$.

\begin{prop}\label{prop3}
Let $(M,g,J)$ be Hermitian, $Hol(\nabla^a) = \rho (Sp(n)U(1))$ and $\nabla^a T^a =0$. Then $T^a = \lambda T_0 + \bar{\lambda} \overline{T}_0$, where $\lambda \in \mathbb{C}$ is a constant and $T_0$ is the $\nabla^a$-parallel tensor field defined by the $\rho (Sp(n)U(1))$-invariant tensor $T_0$ from Section~\ref{sec1}.
\end{prop}

\noindent {\it Proof:}
$\nabla^a T^a =0$ implies that $T^a$ is invariant with respect to $Hol(\nabla^a) = \rho (Sp(n)U(1))$ and the assertion follows from Corollary~\ref{cor1}. \hfill $\Box$

\vspace{3mm}
Let $(M^{2m+2},g,J)$ belong to the class $\mathcal{G}_1$ and $Hol(\nabla^a) \subset U(m) \times U(1)$. Then we have a $\nabla^a$-parallel, orthogonal and $J$-invariant splitting $TM = \mathcal{H} \oplus \mathcal{V}$, where $\dim_{\mathbb{R}} \mathcal{H} = 2m$, $\dim_{\mathbb{R}} \mathcal{V} = 2$. We can define an orthogonal with respect to $g$ almost complex structure $\hat{J}$ by
$$\hat{J}_{|\mathcal{H}} = J_{|\mathcal{H}}, \quad \hat{J}_{|\mathcal{V}}= -J_{|\mathcal{V}}.$$

\begin{prop}\label{prop4}
Let $(M,g,J)$ be a Hermitian manifold with $Hol(\nabla^a) \subset U(m) \times U(1)$ and $T^a \in  \Lambda ^{2,0} \mathcal{H}^* \otimes \Lambda ^{0,1} \mathcal{V}^* \oplus \Lambda ^{0,2} \mathcal{H}^* \otimes \Lambda ^{1,0} \mathcal{V}^*$. Then $(M,g,\hat{J})$ is nearly K\"ahler and its canonical connection coincides with $\nabla^a$. Conversely, if $(M,g,J)$ is nearly K\"ahler with $Hol(\nabla^a) \subset U(m) \times U(1)$ and $T^a \in  \Lambda ^{2,0} \mathcal{H}^* \otimes \Lambda ^{1,0} \mathcal{V}^* \oplus \Lambda ^{0,2} \mathcal{H}^* \otimes \Lambda ^{0,1} \mathcal{V}^*$, then $(M,g,\hat{J})$ is Hermitian and the unique Hermitian connection with totally skew-symmetric torsion coincides with $\nabla^a$.
\end{prop}

\noindent {\bf Remark:} The condition $T^a \in  \Lambda ^{2,0} \mathcal{H}^* \otimes \Lambda ^{1,0} \mathcal{V}^* \oplus \Lambda ^{0,2} \mathcal{H}^* \otimes \Lambda ^{0,1} \mathcal{V}^*$ is automatically satisfied for a 6-dimensional nearly K\"ahler manifold.

\vspace{3mm}
\noindent {\it Proof:} The definition of $\hat{J}$ implies that $\nabla^a \hat{J} = 0$. Hence $\nabla^a$ is the Hermitian connection with totally skew-symmetric torsion also for $(M,g,\hat{J})$.

If $(M,g,J)$ is Hermitian and $T^a \in  \Lambda ^{2,0} \mathcal{H}^* \otimes \Lambda ^{0,1} \mathcal{V}^* \oplus \Lambda ^{0,2} \mathcal{H}^* \otimes \Lambda ^{1,0} \mathcal{V}^*$, then $T^a$ is a $((3,0) \oplus (0,3))$-form with respect to $\hat{J}$ and therefore $(M,g,\hat{J})$ is nearly K\"ahler.

If $(M,g,J)$ is nearly K\"ahler with $T^a \in  \Lambda ^{2,0} \mathcal{H}^* \otimes \Lambda ^{1,0} \mathcal{V}^* \oplus \Lambda ^{0,2} \mathcal{H}^* \otimes \Lambda ^{0,1} \mathcal{V}^*$, then $T^a$ is a $((2,1) \oplus (1,2))$-form with respect to $\hat{J}$ and therefore $(M,g,\hat{J})$ is Hermitian. \hfill $\Box$

\vspace{3mm}
Since the torsion of the canonical connection of a nearly K\"ahler manifold is parallel \cite{K}, we get

\begin{cor}\label{cor2}
If $(M,g,J)$ is a Hermitian manifold with $Hol(\nabla^a) \subset U(m) \times U(1)$ and $T^a \in  \Lambda ^{2,0} \mathcal{H}^* \otimes \Lambda ^{0,1} \mathcal{V}^* \oplus \Lambda ^{0,2} \mathcal{H}^* \otimes \Lambda ^{1,0} \mathcal{V}^*$, then $\nabla^a T^a =0$.
\end{cor}

\noindent {\bf Remark:} The condition $T^a \in  \Lambda ^{2,0} \mathcal{H}^* \otimes \Lambda ^{0,1} \mathcal{V}^* \oplus \Lambda ^{0,2} \mathcal{H}^* \otimes \Lambda ^{1,0} \mathcal{V}^*$ implies that $\Omega$ is co-closed, i.e., $(M,g,J)$ is semi-K\"ahler (or balanced).

\vspace{3mm}
We summarize now some simple facts about Riemannian manifolds $(M,g)$ with metric connection $\nabla^a = \nabla +\frac{1}{2} T^a$ with totally skew-symmetric torsion $T^a$,  such that $\nabla^a T^a =0$. In this case the first Bianchi identity is
\begin{equation}\label{201}
\mathfrak{S}_{X,Y,Z} R^a(X,Y,Z,W) = \sigma_{T^a} (X,Y,Z,W),
\end{equation}
where $\mathfrak{S}_{X,Y,Z}$ denotes a cyclic sum with respect to $X$, $Y$, $Z$ and $\sigma_{T^a} \in \Lambda^4T^* M$ is defined by
\begin{equation}\label{202}
\sigma_{T^a} (X,Y,Z,W) = \mathfrak{S}_{X,Y,Z} \, g(T^a(X,Y),T^a(Z,W)).
\end{equation}
(For the curvature tensors we use the following convention:
$$R(X,Y) = [\nabla_X,\nabla_Y] - \nabla_{[X,Y]}, \quad R(X,Y,Z,W) = g(R(X,Y)Z,W).)$$
The Bianchi identity \pref{201} implies $R^a(X,Y,Z,W) = R^a(Z,W,X,Y)$.

Since $\nabla^a T^a =0$, we have $R^a(U,V)(T^a) = 0$, i.e.,
\begin{equation}\label{203}
\mathfrak{S}_{X,Y,Z} R^a(U,V,X,T^a(Y,Z)) = 0.
\end{equation}
The explicit relation between $R^a$ and the curvature $R$ of the Levi-Civita connection is given by
\begin{eqnarray}\label{204}
& & R^a(X,Y,Z,W) = R(X,Y,Z,W) + \frac{1}{2}g(T^a(X,Y),T^a(Z,W)) \\
& & \quad 
+ \frac{1}{4}g(T^a(Y,Z),T^a(X,W)) - \frac{1}{4}g(T^a(X,Z),T^a(Y,W)). \nonumber
\end{eqnarray}
Therefore we get the following relations between the Ricci tensors and the scalar curvatures of $\nabla^a$ and $\nabla$:
\begin{equation}\label{205}
Ric^a = Ric - \frac{1}{4}r^a, \qquad 
s^a = s - \frac{|T^a|^2}{4},
\end{equation}
where $r^a(X,Y) = g(T^a(X,\cdot ), T^a(Y,\cdot ))$.

\vspace{3mm}
From now on we shall consider almost Hermitian manifolds $(M,g,J)$ belonging to $\mathcal{G}_1$ such that $T^a$ is non-degenerate, i.e., $T^a(X,\cdot ) \not = 0$ for each $X \not = 0$. Notice that for Hermitian manifolds this condition is weaker than the requirement $\nabla_X J \not = 0$ for each $X \not = 0$. For nearly K\"ahler manifolds the two conditions are equivalent and the manifolds satisfying them are called {\it strict} nearly K\"ahler \cite{Gr}. 

\begin{prop}\label{prop5}
Let $(M,g,J)$ be a Hermitian manifold such that $Hol(\nabla^a) \subset U(m) \times U(1)$, $T^a \in  \Lambda ^{2,0} \mathcal{H}^* \otimes \Lambda ^{0,1} \mathcal{V}^* \oplus \Lambda ^{0,2} \mathcal{H}^* \otimes \Lambda ^{1,0} \mathcal{V}^*$ and $T^a$ is non-degenerate. Then $m=2n$, $Hol(\nabla^a) \subset \rho (Sp(n)U(1))$ and $T^a = \lambda T_0 + \bar{\lambda} \overline{T}_0$, where $\lambda \in \mathbb{C} \backslash \{0\}$ is a constant.
\end{prop}

\noindent {\it Proof:}
From Corollary~\ref{cor2} we know that $\nabla^a T^a =0$. Since $\nabla^a J =0$, we have
$$R^a(JX,JY,Z,W) = R^a(X,Y,JZ,JW) = R^a(X,Y,Z,W), \quad Ric^a (JX,JY) = Ric^a (X,Y).$$
Let $\{e_\alpha \}$ be a basis of $T^{1,0} M$. Then the Bianchi identity \pref{201} becomes
\begin{equation}\label{207}
{R^a}_{\alpha \bar{\beta} \gamma \bar{\delta}} - {R^a}_{\gamma \bar{\beta} \alpha \bar{\delta}} = (\sigma_{T^a})_{\alpha \bar{\beta} \gamma \bar{\delta}}.
\end{equation}
Contracting \pref{203} with respect to $V$ and $X$, we get
\begin{equation}\label{208}
{Ric^a}_A^S {T^a}_{BCS} + g^{PQ}g^{RS}({R^a}_{APBR}{T^a}_{CQS} + {R^a}_{APCR}{T^a}_{QBS}) = 0.
\end{equation}
Let $A=\alpha$, $B=\bar{\beta}$, $C=\bar{\gamma}$. Since $T^a$ is a $((2,1) \oplus (1,2))$-form, we obtain
\begin{equation}\label{209}
{Ric^a}_\alpha^\sigma {T^a}_{\bar{\beta} \bar{\gamma} \sigma} + g^{\bar{\delta} \tau}g^{\epsilon \bar{\sigma}}(-{R^a}_{\alpha \bar{\delta} \bar{\beta} \epsilon}{T^a}_{\tau \bar{\gamma} \bar{\sigma}} + {R^a}_{\alpha \bar{\delta} \bar{\gamma} \epsilon}{T^a}_{\tau \bar{\beta} \bar{\sigma}}) = 0.
\end{equation}
In our case $T^a = \omega \wedge \bar{e}^{m+1} + \bar{\omega} \wedge e^{m+1}$, where $\omega \in \Lambda ^{2,0} \mathcal{H}^*$ and $e^{m+1} \in \Lambda ^{1,0} \mathcal{V}^*$, $|e^{m+1}|=1$. Since $T^a$ is non-degenerate, $\omega$ must be a non-degenerate 2-form on $\mathcal{H}^{1,0}$. Hence $\mathcal{H}^{1,0}$ is even dimensional, i.e., $m=2n$. We can take an orthonormal basis $e_1,\dots ,e_{2n+1}$ of $T^{1,0} M$ so that $e_1,\dots ,e_{2n} \in \mathcal{H}^{1,0}$, $e_{2n+1} \in \mathcal{V}^{1,0}$ and $\omega = \sum _{k=1}^n \lambda_k e^{2k-1} \wedge e^{2k}$, where $\lambda_k >0$, $k=1,\dots ,n$.

Let $\alpha = 2n+1$, $\beta = 2k-1$, $\gamma = 2k$ in \pref{209}. Then we get
\begin{equation}\label{210}
\lambda_k {Ric^a}_{2n+1}^{2n+1} -\lambda_k {R^a}_{2n+1 \overline{2n+1} 2k-1 \overline{2k-1}} - \lambda_k {R^a}_{2n+1 \overline{2n+1} 2k \overline{2k}} = 0.
\end{equation}
From \pref{207}
$${R^a}_{2n+1 \overline{2n+1} 2k-1 \overline{2k-1}} = {R^a}_{2k-1 \overline{2n+1} 2n+1 \overline{2k-1}} + (\sigma_{T^a})_{2n+1 \overline{2n+1} 2k-1 \overline{2k-1}}.$$
The splitting $TM = \mathcal{H} \oplus \mathcal{V}$ is $\nabla^a$-parallel and therefore is preserved by $R^a(X,Y)$. Hence ${R^a}_{2k-1 \overline{2n+1} 2n+1 \overline{2k-1}} = 0$. Thus
$${R^a}_{2n+1 \overline{2n+1} 2k-1 \overline{2k-1}} = (\sigma_{T^a})_{2n+1 \overline{2n+1} 2k-1 \overline{2k-1}} = {\lambda_k}^2.$$
Similarly, ${R^a}_{2n+1 \overline{2n+1} 2k \overline{2k}} = {\lambda_k}^2$. Hence \pref{210} yields
$${Ric^a}_{2n+1}^{2n+1} = 2{\lambda_k}^2, \quad k=1,\dots ,n.$$
This means that $\lambda_1 = \dots = \lambda_n =: \lambda >0$ and
$$T^a = \lambda \sum _{k=1}^n e^{2k-1} \wedge e^{2k} \wedge \bar{e}^{2n+1} + \lambda \sum _{k=1}^n \bar{e}^{2k-1} \wedge \bar{e}^{2k} \wedge e^{2n+1}.$$
Now Proposition~\ref{prop2} implies that $Hol(\nabla^a)$ is conjugate to a subgroup of $\rho (Sp(n)U(1))$. \hfill $\Box$

\vspace{3mm}
The computations in the above proof yield the following.

\begin{cor}\label{cor3}
Let $(M,g,J)$ be as in Proposition~\ref{prop5}. Then $Ric^a$, $s^a$, $Ric$, $s$ are given by \pref{211} and \pref{212} below. In particular, $Ric^a$ and $Ric$ are positive definite and $\nabla^a$-parallel and $(M,g)$ is Einstein only if $n=1$.
\end{cor}

\noindent {\it Proof:}
From \pref{208} with $A=\alpha$, $B=\bar{\beta}$, $C=\gamma$ and \pref{209} in a similar way as in the proof of Proposition~\ref{prop5} we get
$${Ric^a}_\alpha^\sigma = \left\{
\begin{array}{ll}
(n+1)\lambda^2, & \quad \alpha = \sigma < 2n+1, \\
2\lambda^2, & \quad \alpha = \sigma = 2n+1, \\
0, & \quad \alpha \not = \sigma .
\end{array}
\right.$$
This and \pref{101} yield \pref{211}. \pref{212} follows from \pref{211} and \pref{205}. \hfill $\Box$

\vspace{3mm}
\noindent {\bf Remark:} Let $(M,g,J)$ be a strict nearly K\"ahler manifold. Then $T^a$ is a $((3,0) \oplus (0,3))$-form. Taking $A=\alpha$, $B=\beta$, $C=\gamma$ in \pref{208} we get
$${Ric^a}_\alpha^\sigma {T^a}_{\beta \gamma \sigma} + g^{\bar{\delta} \tau}g^{\bar{\epsilon} \sigma}({R^a}_{\alpha \bar{\delta} \beta \bar{\epsilon}}{T^a}_{\gamma \tau \sigma} - {R^a}_{\alpha \bar{\delta} \gamma \bar{\epsilon}}{T^a}_{\beta \tau \sigma}) = 0.$$
This together with \pref{207} yields
$${Ric^a}_\alpha^\sigma {T^a}_{\beta \gamma \sigma} + g^{\bar{\delta} \tau}g^{\bar{\epsilon} \sigma}((\sigma_{T^a})_{\alpha \bar{\delta} \beta \bar{\epsilon}}{T^a}_{\gamma \tau \sigma} - (\sigma_{T^a})_{\alpha \bar{\delta} \gamma \bar{\epsilon}}{T^a}_{\beta \tau \sigma}) = 0,$$
which in global notations can be written as
\begin{equation}\label{214}
T^a(Y,Z,Ric^a(X)) = \frac{1}{2}T^a(X,r^a(Y),Z) + \frac{1}{2}T^a(X,Y,r^a(Z)).
\end{equation}
Hence $T^a(Y,Z,(\nabla^a_W Ric^a)(X)) = 0$ for each $Y,Z$. Since $T^a$ is non-degenerate, this implies $(\nabla^a_W Ric^a)(X) = 0$. Thus $\nabla^a Ric^a = 0$ and by \pref{205} also $\nabla^a Ric = 0$.

This proof is essentially due to Kirichenko \cite{K} but he concludes wrongly from \pref{214} that the manifold is Einstein. The correct formulation can be found in \cite{N}. From \pref{214} one can also see that $Ric^a$ and $r^a$ commute. Hence, there is an orthonormal basis $\{ e_\alpha \}$ of $T^{1,0} M$ consisting of eigenvectors for both $Ric^a$ and $r^a$. If the corresponding eigenvalues are $\mu_\alpha$ and $\nu_\alpha$, then \pref{214} shows that $\mu_\alpha = \frac{1}{2}(\nu_\beta + \nu_\gamma)$ whenever ${T^a}_{\alpha \beta \gamma} \not =0$. Thus, since $r^a$ is positive definite, we get that $Ric^a$ (and by \pref{205} also $Ric$) is positive definite. This has also been proved in \cite{N}.

\vspace{3mm}
From Proposition~\ref{prop4} and Proposition~\ref{prop5} we get
\begin{cor}\label{cor4}
Let $(M,g,J)$ be a strict nearly K\"ahler manifold with $Hol(\nabla^a) \subset U(m) \times U(1)$ and $T^a \in  \Lambda ^{2,0} \mathcal{H}^* \otimes \Lambda ^{1,0} \mathcal{V}^* \oplus \Lambda ^{0,2} \mathcal{H}^* \otimes \Lambda ^{0,1} \mathcal{V}^*$. Then $m=2n$ and $Hol(\nabla^a) \subset \rho_2 (Sp(n)U(1))$.
\end{cor}

\section{The curvature}\label{sec5}

In section~\ref{sec2} we saw that if $(M,g,J)$ is a Hermitian manifold with $Hol(\nabla^a) = \rho (Sp(n)U(1))$ and $\nabla^a T^a =0$, then $T^a = \lambda T_0 + \bar{\lambda} \overline{T}_0$ and the curvature tensor $R^a$ has the following properties:

$R^a (X,Y) \in \rho (sp(n) \oplus u(1))$,

$R^a$ is symmetric with respect to first and second pair of arguments,

$bR^a = \sigma_{T^a} $, \\
where $b:(T^*)^{\otimes 4} \longrightarrow (T^*)^{\otimes 4}$ is
$$(bR) (X,Y,Z,W) = \mathfrak{S}_{X,Y,Z} R(X,Y,Z,W).$$
Let $\mathfrak{R}^a$ denote the space of algebraic tensors with these properties, i.e.,
$$\mathfrak{R}^a = S^2 (\rho (sp(n) \oplus u(1))) \cap b^{-1} (\mathbb{R} \sigma_{T_0 + \overline{T}_0}).$$

\begin{prop}\label{prop7}
The complexification of $\mathfrak{R}^a$ is isomorphic to $\mathbb{C} \oplus S^4 E$ as a $\rho (Sp(n)U(1))$-representation.
\end{prop}

\noindent {\it Proof:}
$\mathbb{R} \sigma_{T_0 + \overline{T}_0}$ and $u(1) \cong \mathbb{R}$ are trivial representations of $\rho (Sp(n)U(1))$ and the complexification of $sp(n)$ is isomorphic to $S^2 E$. Hence the complexification of $\mathfrak{R}^a$ is
$$\mathfrak{R}^a \otimes \mathbb{C} \cong S^2 (S^2 E \oplus \mathbb{C}) \cap b^{-1} (\mathbb{C} \sigma_{T_0 + \overline{T}_0}).$$
We have
$$S^2 (S^2 E \oplus \mathbb{C}) \cong S^2 (S^2 E) \oplus S^2 E \oplus \mathbb{C}$$
and
$$S^2 (S^2 E)\cong S^4 E \oplus L \oplus \Lambda^2_0 E \oplus \mathbb{C},$$
where $L$ is the representation of $Sp(n)$ with highest weight $(2,2,0,\dots ,0)$ (if $n=1$, $L$ and $\Lambda^2_0 E$ are 0). Thus the decomposition of $S^2 (S^2 E \oplus \mathbb{C})$ into irreducible $\rho (Sp(n)U(1))$-representations is
$$S^2 (S^2 E \oplus \mathbb{C}) \cong S^4 E \oplus L \oplus \Lambda^2_0 E \oplus S^2 E \oplus \mathbb{C} \oplus \mathbb{C}.$$
Taking particular representatives of these spaces and using Schur's Lemma one sees that $S^4 E \subset \ker b$, that $b$ is non-zero on $L$, $\Lambda^2_0 E$, $S^2 E$ and therefore they are not contained in $\mathfrak{R}^a \otimes \mathbb{C}$, and that $b_{|\mathbb{C} \oplus \mathbb{C}}$ is injective and $b(\mathbb{C} \oplus \mathbb{C}) \supset \mathbb{C} \sigma_{T_0 + \overline{T}_0}$. Thus $\mathfrak{R}^a \otimes \mathbb{C} \cong \mathbb{C} \oplus S^4 E$. \hfill $\Box$

\begin{cor}\label{cor6}
As a real representation of $\rho (Sp(n)U(1))$, $\mathfrak{R}^a \cong \mathbb{R} \oplus r(S^4 E)$, where $r(S^4 E)$ is the real representation underlying $S^4 E$. If $R^a \in \mathfrak{R}^a$ and the corresponding torsion is $T^a = \lambda T_0 + \bar{\lambda} \overline{T}_0$, then with respect to the above isomorphism
\begin{equation}\label{307}
R^a = \frac{|T^a|^2}{48n} R^a_0 + R_{hyper},
\end{equation}
where $R_{hyper}$ has the properties of an algebraic hyper-K\"ahler curvature tensor on $\mathcal{H}$ and
\begin{eqnarray}\label{308}
& & R^a_0 (X,Y,Z,W) =  \\
& & \qquad \qquad \sum_{L \in \{ \pmb{1},I,J,K \}} (g_{|\mathcal{H}}(LY,Z)g_{|\mathcal{H}}(LX,W) - g_{|\mathcal{H}}(LX,Z)g_{|\mathcal{H}}(LY,W)) \nonumber \\
& & \qquad \qquad - 2(g_{|\mathcal{H}}(JX,Y)g_{|\mathcal{H}}(JZ,W) + 2g_{|\mathcal{H}}(JX,Y)g_{|\mathcal{V}}(JZ,W) \nonumber \\
& & \qquad \qquad \qquad \qquad + 2g_{|\mathcal{V}}(JX,Y)g_{|\mathcal{H}}(JZ,W) + 4g_{|\mathcal{V}}(JX,Y)g_{|\mathcal{V}}(JZ,W)). \nonumber
\end{eqnarray}
The Ricci tensor $Ric^a$ and the scalar curvature $s^a$ of $R^a$ are 
\begin{equation}\label{211}
Ric^a = \frac{|T^a|^2}{12n} ((n+1)g_{|\mathcal{H}} + 2g_{|\mathcal{V}}), \quad s^a = \frac{(n^2+n+1)|T^a|^2}{3n}.
\end{equation}
\end{cor}

\noindent {\it Proof:}
One needs only to check that the preimage in $\mathbb{C} \oplus \mathbb{C}$ of $\sigma_{T^a}$ with respect to $b$ is $\frac{|T^a|^2}{48n} R^a_0$ (recall that $|T^a|^2 = 12n |\lambda|^2$). The formulae for $Ric^a$ and $s^a$ follow from the explicit form of $R^a_0$ (or, alternatively, from Corollary~\ref{cor3}).  \hfill $\Box$

\vspace{3mm}
From \pref{204} we obtain

\begin{cor}\label{cor7}
As a $\rho (Sp(n)U(1))$-representation, the space $\mathfrak{R}$ of algebraic tensors with the properties of the Riemannian curvature tensor of a Hermitian manifold with $Hol(\nabla^a) = \rho (Sp(n)U(1))$ and $\nabla^a T^a =0$ is $\mathfrak{R} \cong \mathbb{R} \oplus r(S^4 E)$. With respect to this isomorphism, for $R \in \mathfrak{R}$, we have
\begin{equation}\label{305}
R = \frac{|T^a|^2}{48n} R_0 + R_{hyper},
\end{equation}
where $R_{hyper}$ has the properties of an algebraic hyper-K\"ahler curvature tensor on $\mathcal{H}$  and
\begin{eqnarray}\label{306}
& & R_0 (X,Y,Z,W) =  R_{\mathbb{H}P^n}^{g_{|\mathcal{H}}} (X,Y,Z,W)     \\
& & \qquad -\frac{1}{2} \sum_{L \in \{ I,K \}} ( g_{|\mathcal{H}}(LY,Z)g_{|\mathcal{H}}(LX,W) - g_{|\mathcal{H}}(LX,Z)g_{|\mathcal{H}}(LY,W)  \nonumber \\
& & \qquad \qquad \qquad \qquad - 2g_{|\mathcal{H}}(LX,Y)g_{|\mathcal{H}}(LZ,W) ) \nonumber \\
& & \qquad + \frac{1}{2}(g_{|\mathcal{H}}(Y,Z)g_{|\mathcal{V}}(X,W) + g_{|\mathcal{V}}(Y,Z)g_{|\mathcal{H}}(X,W) \nonumber \\
& & \qquad \qquad \qquad - g_{|\mathcal{H}}(X,Z)g_{|\mathcal{V}}(Y,W)- g_{|\mathcal{V}}(X,Z)g_{|\mathcal{H}}(Y,W)) \nonumber \\
& & \qquad + \frac{3}{2}(g_{|\mathcal{H}}(JY,Z)g_{|\mathcal{V}}(JX,W) + g_{|\mathcal{V}}(JY,Z)g_{|\mathcal{H}}(JX,W) \nonumber \\
& & \qquad \qquad \qquad - g_{|\mathcal{H}}(JX,Z)g_{|\mathcal{V}}(JY,W)- g_{|\mathcal{V}}(JX,Z)g_{|\mathcal{H}}(JY,W) \nonumber \\
& & \qquad \qquad \qquad - 2g_{|\mathcal{H}}(JX,Y)g_{|\mathcal{V}}(JZ,W) - 2g_{|\mathcal{V}}(JX,Y)g_{|\mathcal{H}}(JZ,W)) \nonumber \\
& & \qquad -8g_{|\mathcal{V}}(JX,Y)g_{|\mathcal{V}}(JZ,W) \nonumber
\end{eqnarray}
($R_{\mathbb{H}P^n}^{g_{|\mathcal{H}}}$ is given by \pref{300} with $g'$ replaced by $g_{|\mathcal{H}}$ and $I'$, $J'$, $K'$ replaced by $I$, $J$, $K$). The Ricci tensor $Ric$ and the scalar curvature $s$ of $R$ are 
\begin{equation}\label{212}
Ric = \frac{|T^a|^2}{24n} ((2n+3)g_{|\mathcal{H}} + (n+4)g_{|\mathcal{V}}), \quad s = \frac{(4n^2+7n+4)|T^a|^2}{12n}.
\end{equation}
\end{cor}

Proposition~\ref{prop5} yields

\begin{cor}\label{cor8}
If $(M,g,J)$ is as in Proposition~\ref{prop5}, then $R^a$ and $R$ are as in Corollary~\ref{cor6} and Corollary~\ref{cor7} and $\nabla^a R^a_0 = 0$, $\nabla^a R_0 = 0$.
\end{cor}

\noindent {\bf Remark:} The results in this section have obvious counterparts for strict nearly K\"ahler manifolds with $Hol(\nabla^a) \subset U(m) \times U(1)$ and $T^a \in  \Lambda ^{2,0} \mathcal{H}^* \otimes \Lambda ^{1,0} \mathcal{V}^* \oplus \Lambda ^{0,2} \mathcal{H}^* \otimes \Lambda ^{0,1} \mathcal{V}^*$ (and therfore with $Hol(\nabla^a) \subset \rho_2 (Sp(n)U(1))$ by Corollary~\ref{cor4}). One only needs to change the signs of the summands of the type $g_{|\mathcal{H}} (J \cdot ,\cdot )  g_{|\mathcal{V}} (J \cdot ,\cdot )$ in the formulae for $R^a_0$ and $R_0$.

\vspace{3mm}
\noindent {\bf Remark:} The method of the proof of Proposition~\ref{prop7} can be used to prove the already mentioned fact that the canonical connection of a nearly K\"ahler manifold has parallel torsion. In general, the first Bianchi identity has the form
$${\mathfrak{S}}_{X,Y,Z} R^a(X,Y,Z,W) =  {\mathfrak{S}}_{X,Y,Z}(\nabla^a_X T^a (Y,Z,W) + g(T^a(X,Y),T^a(Z,W)).$$
Since $Hol(\nabla^a) \subset U(n)$, $R^a \in \Lambda^2 T^* \otimes u(n)$. Furthermore, $T^a$ is a $((3,0) \oplus (0,3))$-form and therefore $\nabla^a T^a \in T^* \otimes (\Lambda^{3,0} T^* \oplus \Lambda^{0,3} T^*)$.  Let $\tau^a (X,Y,Z,W) = g(T^a(X,Y),T^a(Z,W))$. Then
$$\tau^a (X,Y,Z,W) = g((T^a)^{3,0}(X,Y),(T^a)^{0,3}(Z,W)) + g((T^a)^{0,3}(X,Y),(T^a)^{3,0}(Z,W)).$$
Hence, complexifying all spaces, we obtain
$$b(R^a-\nabla^a T^a - \tau^a) =0$$
with $R^a \in \Lambda^2 (T^\mathbb{C})^* \otimes u(n)^\mathbb{C}$, $\nabla^a T^a \in ((T^{1,0})^* \oplus (T^{0,1})^*) \otimes (\Lambda^{3,0} T^* \oplus \Lambda^{0,3} T^*)$, $\tau^a \in \Lambda^{2,0 }T^* \otimes \Lambda^{0,2} T^*$. Decomposing these spaces into irreducible $U(n)$-representations and taking particular representatives and using Schur's Lemma to determine the rank of $b$ on the different components, we obtain
\begin{eqnarray}\label{3001}
& & ((T^{1,0})^* \oplus (T^{0,1})^*) \otimes (\Lambda^{3,0} T^* \oplus \Lambda^{0,3} T^*) = \Lambda^{4,0} T^* \oplus \Lambda^{0,4} T^* \\
& & \qquad \oplus V(2,1,1,0,\dots ,0) \otimes F(4) \oplus V(0,\dots ,0,-1,-1,-2) \otimes F(-4) \nonumber \\
& & \qquad \oplus \Lambda^{2,0} T^* \oplus \Lambda^{0,2} T^* \nonumber \\
& & \qquad \oplus V(1,1,1,0,\dots ,0,-1) \otimes F(2) \oplus V(1,0,\dots ,0,-1,-1,-1) \otimes F(-2) \nonumber
\end{eqnarray}
($V(\alpha)$ is the irreducible representation of $SU(n)$ with highest weight $\alpha$, $F(k)$ is the representation of $U(1)$ with weight $k$) and
\begin{itemize}
\item $b$ is injective on the first four components in \pref{3001} and they are not contained in $\Lambda^2 (T^\mathbb{C})^* \otimes u(n)^\mathbb{C}$ and $\Lambda^{2,0} T^* \otimes \Lambda^{0,2} T^*$. Therefore $\nabla^a T^a$ has no components in them.
\item $\Lambda^{2,0} T^*$ is contained twice in $\Lambda^2 (T^\mathbb{C})^* \otimes u(n)^\mathbb{C}$, but $b_{|\Lambda^{2,0} T^* \oplus \Lambda^{2,0} T^* \oplus \Lambda^{2,0} T^*}$ is injective. Hence $\nabla^a T^a$ has no component in $\Lambda^{2,0} T^*$. In a similar way $\nabla^a T^a$ has no component in $\Lambda^{0,2} T^*$.
\item $V(1,1,1,0,\dots ,0,-1) \otimes F(2)$ \hfill is \hfill contained \hfill once \hfill in \hfill $\Lambda^2 (T^\mathbb{C})^* \otimes u(n)^\mathbb{C}$, \hfill but \\
$b_{|V(1,1,1,0,\dots ,0,-1) \otimes F(2) \oplus V(1,1,1,0,\dots ,0,-1) \otimes F(2)}$ is injective. Thus $\nabla^a T^a$ has no component in $V(1,1,1,0,\dots ,0,-1) \otimes F(2)$. The same is true for $V(1,0,\dots ,0,-1,-1,-1) \otimes F(-2)$.
\end{itemize}

Hence $\nabla^a T^a = 0$.

\section{Examples: the twistor spaces}\label{sec3}

Recall that a quaternionic K\"ahler manifold is a $4n$-dimensional Riemannian manifold $(M',g')$ whose holonomy is contained in $Sp(n)Sp(1)$ if $n>1$ or which is self-dual and Einstein if $n=1$. If $n>1$ an equivalent definition is to require the existence of a subbundle $Q' \subset End(TM')$ of rank 3 which is locally trivialized by three orthogonal almost complex structures $I'$, $J'$, $K'$ satisfying the quaternionic identities. Such a bundle exists also if $n=1$. In this case we choose $Q' = \Lambda^2_- M'$ (if we would like to choose the other possibility $Q' = \Lambda^2_+ M'$, we have to replace "self-dual" by "anti-self-dual" in the definition of quaternionic K\"ahler manifold).

Every quaternionic K\"ahler manifold is Einstein and its curvature has the form
$$R' = \frac{s'}{16n(n+2)} R_{\mathbb{H}P^n} + R'_{hyper},$$
where $s'$ is the (constant) scalar curvature, $R_{\mathbb{H}P^n}$ is the (parallel) curvature tensor of $\mathbb{H}P^n$,
\begin{eqnarray}\label{300}
R_{\mathbb{H}P^n} (X,Y,Z,W) = g'(Y,Z)g'(X,W) - g'(X,Z)g'(Y,W) \\
+ \sum_{L \in \{I',J',K'\}} (g'(LY,Z)g'(LX,W) - g'(LX,Z)g'(LY,W) - 2g'(LX,Y)g'(LZ,W)), \nonumber
\end{eqnarray}
and $R'_{hyper}$ has the symmetries of a hyper-K\"ahler curvature tensor. If $n=1$,  $R_{\mathbb{H}P^n}$ is the curvature of $S^4$ with the metric with sectional curvature 4 and $R'_{hyper} = W_+$ (the positive Weyl tensor).

The twistor space $\mathcal{Z}$ of a quaternionic K\"ahler manifold $M'$ is the $S^2$-bundle over $M'$ whose fibre at $p \in M'$ is $\mathcal{Z}_p = \{ z \in Q'_p:z^2 = -\pmb{1} \}$. A local trivialization $\psi = (\pi,\varphi)$ of $\mathcal{Z}$ is defined by a local frame $I'$, $J'$, $K'$ of $Q'$, which  satisfies the quaternionic identities, as follows:
if $z \in \mathcal{Z}_p$, $z=aI'+bJ'+cK'$, then $\varphi(z) = (a,b,c) \in S^2 \subset \mathbb{R}^3$.

The Levi-Civita connection defines a horizontal distribution $\mathcal{H}$ on $\mathcal{Z}$. Let $\mathcal{V}$ be the vertical distribution (tangent to the fibres). Two almost complex structures $J_1$ and $J_2$ and a one-parameter family of Riemannian metrics $h_t$, $t>0$, are defined on $\mathcal{Z}$ in the following way: at $z \in \mathcal{Z}$

${J_1}_{|\mathcal{H}_z} = {J_2}_{|\mathcal{H}_z}$ is the complex structure corresponding to $z$ under the isomorphism of $\mathcal{H}_z$ and $T_{\pi (z)} M'$ given by the projection $\pi: \mathcal{Z} \longrightarrow M'$,

${J_1}_{|\mathcal{V}_z} = -{J_2}_{|\mathcal{V}_z}$ corresponds to the standard complex structure of $S^2$ via $\psi$,

${h_t}_{|\mathcal{H}_z}$ corresponds to $g'_{\pi(z)}$ via $\pi$,

${h_t}_{|\mathcal{V}_z}$ corresponds via $\psi$ to the metric with sectional curvature $\frac{1}{nt}$ on $S^2$, 

$\mathcal{H}$ and $\mathcal{V}$ are orthogonal with respect to $h_t$.

\vspace{2mm}
Every two frames of $Q'$ satisfying the quaternionic identities are related by an $SO(3)$-matrix and therefore the definition of $J_1$, $J_2$ and $h_t$ is independent of the choice of $\psi$.

It is well known that $J_1$ is integrable and $J_2$ is not, that the Riemannian submersions $\pi: (\mathcal{Z},h_t) \longrightarrow (M',g')$ have totally geodesic fibres and that $J_1$ and $J_2$ are orthogonal with respect to $h_t$.

The Hermitian structures $(h_t,J_1)$ are semi-K\"ahler for each $t$ (see \cite{M,AGI}). If the scalar curvature $s'$ of $M'$ is positive, there exist two especially interesting values of the parameter $t$:
$(\mathcal{Z},h_t,J_1)$ is K\"ahler iff $t=t_0:=\frac{4(n+2)}{s'}$ and
$(\mathcal{Z},h_t,J_2)$ is nearly K\"ahler iff $t=t_1:=\frac{2(n+2)}{s'}$
(see \cite{FK,S,M,AGI}).

Let us consider $\nabla^{a,t}$, the Hermitian connection with totally skew-symmetric torsion  of $(\mathcal{Z},h_t,J_1)$. An immediate consequence of the results in \cite{M,AGI} is that its torsion $T^{a,t}$ is given by
\begin{equation}\label{301}
T^{a,t} = \frac{1}{\sqrt{2nt}}\left( 2 - \frac{s't}{2(n+2)} \right) (\omega \wedge \bar{\alpha} + \bar{\omega} \wedge \alpha),
\end{equation}
where $\omega \in \Lambda ^{2,0} \mathcal{H}^*$ and  $\alpha \in \Lambda ^{1,0} \mathcal{V}^*$ are defined as follows: \\
Using the trivialization $\psi$ we can consider the vertical vectors at $z \in \mathcal{Z}$ as elements of $Q'_{\pi (z)}$ by identifying $(a,b,c) \in \mathbb{R}^3$ with $aI'+bJ'+cK' \in Q'_{\pi (z)}$. Thus, if $U \in \mathcal{V}_z$, we can define a 2-form $\Omega_U$ on $T_{\pi (z)} M'$ by $\Omega_U (X,Y) = g'(U(X),Y)$. We fix $U \in \mathcal{V}^{1,0}$ with $|U|=1$ and take $\alpha \in \Lambda ^{1,0} \mathcal{V}^*$ to be the dual form of $U$ and $\omega = \frac{\sqrt{2nt}}{2} \pi^* \Omega_U$.

We have
\begin{equation}\label{302}
|T^{a,t}|^2 = \frac{6}{t} \left( 2 - \frac{s't}{2(n+2)} \right) ^2.
\end{equation}
Now, using \cite{M,AGI}, one can prove that $\nabla^{a,t} T^{a,t} =0$ iff $s'>0$ and $t=t_0$ (in this case $T^{a,t_0} =0$) or $t=t_1$. Furthermore, the splitting $T\mathcal{Z} = \mathcal{H} \oplus \mathcal{V}$ is $\nabla^{a,t_1}$-parallel (in fact, the splitting $T\mathcal{Z} = \mathcal{H} \oplus \mathcal{V}$ is always parallel with respect to the canonical Hermitian connection of $(h_t,J_2)$, never parallel with respect to the canonical Hermitian connection of $(h_t,J_1)$ and parallel with respect to $\nabla^{a,t}$ only if $s'>0$ and $t=t_1$). Thus, by Proposition~\ref{prop2} (or Proposition~\ref{prop5}) $Hol(\nabla^{a,t_1}) \subset \rho (Sp(n)U(1))$.

\begin{prop}\label{prop6}
Let $(M',g')$ be a quaternionic K\"ahler manifold with positive scalar curvature. Then the connection $\nabla^{a,t_1}$ on the twistor space $(\mathcal{Z},h_{t_1},J_1)$ has parallel torsion and  $Hol(\nabla^{a,t_1}) \subset \rho (Sp(n)U(1))$. If $(M',g')$ is not locally symmetric, then  $Hol(\nabla^{a,t_1}) = \rho (Sp(n)U(1))$.
\end{prop}

\noindent {\it Proof:}
Only the last assertion remains to be proved. 

Since $Hol(\nabla^{a,t_1}) \subset \rho (Sp(n)U(1))$, we have a $\nabla^{a,t_1}$-parallel quaternionic structure on $\mathcal{H}$, defined as in Section~\ref{sec1}. On the other hand, since ${\pi_*}_{|\mathcal{H}_z} : \mathcal{H}_z \longrightarrow T_{\pi (z)} M'$ is an isomorphism, the quaternionic structure on $M'$ also defines a quaternionic structure on $\mathcal{H}$. We are going to show that these two quaternionic structures coincide.

Fix $z \in \mathcal{Z}$ and choose the trivialization $I'$, $J'$, $K'$ of $Q'$ so that $J'_{\pi (z)} =z$. Then the vector $U \in \mathcal{V}^{1,0}_z$ in the definition of $\alpha$ and $\omega$ above can be taken to be $U =\frac{1}{\sqrt{2nt_1}} (K' + iI')$ and therefore $\omega = \frac{1}{2} \pi^* \Omega_{K'+iI'}$. Thus the quaternionic structure $span \{ I_{|\mathcal{H}_z}, J_{|\mathcal{H}_z}, K_{|\mathcal{H}_z} \}$ on $\mathcal{H}_z$, \hfill defined \hfill by \hfill $Hol(\nabla^{a,t_1})$, \hfill is \hfill exactly \hfill the \hfill pull-back \\
$span \{ ({\pi_*}_{|\mathcal{H}_z})^{-1} (I'), ({\pi_*}_{|\mathcal{H}_z})^{-1} (J'), ({\pi_*}_{|\mathcal{H}_z})^{-1} (K') \}$ of the quaternionic structure $Q'$ on $M'$.

In order to simplify the expressions, from now on we write $g$ instead of $h_{t_1}$, $J$ instead of $J_1$, $\nabla^a$ instead of $\nabla^{a,t_1}$ and $T^a$ instead of $T^{a,t_1}$.

Since $t_1 = \frac{2(n+2)}{s'}$, \pref{301} and \pref{302} become
\begin{equation}\label{303}
T^a = \frac{1}{\sqrt{2nt_1}}(\omega \wedge \bar{\alpha} + \bar{\omega} \wedge \alpha),
\end{equation}
\begin{equation}\label{304}
|T^a|^2 = \frac{3s'}{n+2}.
\end{equation}
Taking into account \pref{304} and that the above defined quaternionic structures on $\mathcal{H}$ coincide, the explicit formulae for the curvature tensor $R$ of the Levi-Civita connection of $(\mathcal{Z},h_t)$ in \cite{DM,AGI}, applied for $t=t_1$, show that $R$ is given by \pref{305},
where $R_{hyper}$ is the horizontal lift of $R'_{hyper}$ (and also  $R_{\mathbb{H}P^n}^{g_{|\mathcal{H}}}$ is the horizontal lift of $R_{\mathbb{H}P^n}$). Now \pref{303} and \pref{204} imply that the curvature of $\nabla^a$ is given by \pref{307}, $R_{hyper}$ being the same as above.

The Lie algebra $hol(\nabla^a)$ of $Hol(\nabla^a)$ is contained in $\rho (sp(n) \oplus u(1))$. On the other hand, $hol(\nabla^a)$ contains the algebra generated by $\{ (\nabla^a)^k_{X_1,\dots ,X_k} R^a(X,Y): k \geq 0 \}$ (in fact,  $\nabla^a$ is real analytic and therefore $hol(\nabla^a)$ is equal to this algebra but we do not need the real analyticity here).

From the definition of $\rho$ we see that $\rho (u(1))$ is spanned by $J_{|\mathcal{H}} + 2J_{|\mathcal{V}}$. If $e_1,\dots ,e_{4n+2}$ is an orthonormal frame of $T\mathcal{Z}$, then \pref{307} and {\pref{308} give
$$J_{|\mathcal{H}} + 2J_{|\mathcal{V}} = -\frac{12n}{(2n+1)|T^a|^2} \sum_{k=1}^{4n+2} R^a (e_k,Je_k).$$
Hence $\rho (u(1)) \subset hol(\nabla^a)$. Therefore to prove that $\rho (sp(n)) \subset hol(\nabla^a)$ it will be enough to show that the algebra generated by
$$A = \{ \mbox{the $\rho (sp(n))$-part of $R^a(X^h,Y^h)$, } (\nabla^a)^k_{X_1^h,\dots ,X_k^h} R^a(X^h,Y^h): k \geq 1 \}$$
contains $\rho (sp(n))$.

The $\rho (sp(n))$-part of $R^a(X^h,Y^h)$ is
$$\frac{|T^a|^2}{48n} \sum_{L \in \{ \pmb{1},I,J,K \}} (g_{|\mathcal{H}}(LY^h,\cdot )LX^h - g_{|\mathcal{H}}(LX^h,\cdot )LY^h) + R_{hyper}(X^h,Y^h)$$
and, because of \pref{304}, this projects on
$$\frac{s'}{16n(n+2)} \sum_{L \in \{ \pmb{1},I',J',K' \}} (g'(LY,\cdot )LX - g'(LX,\cdot )LY) + R'_{hyper}(X,Y),$$
which is exactly the $sp(n)$-part of $R'(X,Y)$.

Since $\pi$ is a Riemannian submersion, $\nabla$ projects on the Levi-Civita connection $\nabla'$ of $g'$, i.e., $h\nabla_{X^h} Y^h = (\nabla'_X Y)^h$. $T^a$ has no components in $\Lambda^3 \mathcal{H}^*$ and therefore $h\nabla^a _{X^h} Y^h = h\nabla_{X^h} Y^h$. Thus $\nabla^a$ also projects on $\nabla'$ (see the Appendix) and since $R_{hyper}$ projects on $R'_{hyper}$, $(\nabla^a)^k R^a = (\nabla^a)^k R_{hyper}$  for $k \geq 1$ projects on $\nabla'^k R'_{hyper}$ by Proposition~\ref{propa1000}.

Hence $A$ projects on
$$B = \{ \mbox{the $sp(n)$-part of $R'(X,Y)$, } \nabla'^k_{X_1,\dots ,X_k} R'_{hyper} (X,Y): k \geq 1 \}.$$

$(M',g')$ is quaternionic K\"ahler with non-zero scalar curvature which is not locally symmetric. Therefore $Hol(\nabla') = Sp(n)Sp(1)$. Since every quaternionic K\"ahler manifold is real analytic,  $hol(\nabla') = sp(n) \oplus sp(1)$ is  generated by $\{ \nabla'^k_{X_1,\dots ,X_k} R' (X,Y): k \geq 0 \}$. But for $k \geq 1$ we have $\nabla'^k_{X_1,\dots ,X_k} R' (X,Y) = \nabla'^k_{X_1,\dots ,X_k} R'_{hyper} (X,Y) \in sp(n)$ and this implies that $sp(n)$ is generated by $B$. Hence  $\rho (sp(n))$ is contained in the algebra generated by $A$. Thus $\rho (sp(n)) \subset hol(\nabla^a)$ and therefore $hol(\nabla^a) = \rho (sp(n) \oplus u(1))$, i.e., $Hol(\nabla^a) = \rho (Sp(n)U(1))$. \hfill $\Box$

\begin{cor}\label{cor5}
Let $(M',g')$ be a quaternionic K\"ahler manifold with positive scalar curvature. Then the canonical connection $\nabla^{a,t_1}$ of the nearly K\"ahler manifold $(\mathcal{Z},h_{t_1},J_2)$ has $Hol(\nabla^{a,t_1}) \subset \rho_2 (Sp(n)U(1))$. If $(M',g')$ is not locally symmetric, then  $Hol(\nabla^{a,t_1}) = \rho_2 (Sp(n)U(1))$.
\end{cor}

\noindent {\bf Remark:} A locally symmetric quaternionic K\"ahler manifold $M'$ has holonomy group $HSp(1)$ where $H \subset Sp(n)$ (see \cite{B} for the list of possible groups $H$). The proof of Proposition~\ref{prop6} shows that in this case the connection $\nabla^{a,t_1}$ on the twistor space $(\mathcal{Z},h_{t_1},J_1)$ has $Hol(\nabla^{a,t_1}) = \rho (HU(1))$. In particular, for $M' = \mathbb{H}P^n$ again $Hol(\nabla^{a,t_1}) = \rho (Sp(n)U(1))$. A similar remark is true for  $(\mathcal{Z},h_{t_1},J_2)$ (replace $\rho$ by $\rho_2$).

\section{The compact case}\label{sec4}

\begin{theo}\label{thm1}
Let $(M,g,J)$ be a complete Hermitian manifold such that $Hol(\nabla^a) \subset U(m) \times U(1)$, $T^a \in  \Lambda ^{2,0} \mathcal{H}^* \otimes \Lambda ^{0,1} \mathcal{V}^* \oplus \Lambda ^{0,2} \mathcal{H}^* \otimes \Lambda ^{1,0} \mathcal{V}^*$ and $T^a$ is non-degenerate. Then $m=2n$ and $(M,g,J)$ is isomorphic to the twistor space $(\mathcal{Z},h_{t_1},J_1)$ of some compact quaternionic K\"ahler manifold with positive scalar curvature.
\end{theo}

\noindent {\it Proof:}
By Proposition~\ref{prop4} $(M,g,\hat{J})$ is a complete strict nearly K\"ahler manifold with  $Hol(\nabla^a) \subset U(m) \times U(1)$ and $T^a \in  \Lambda ^{2,0} \mathcal{H}^* \otimes \Lambda ^{1,0} \mathcal{V}^* \oplus \Lambda ^{0,2} \mathcal{H}^* \otimes \Lambda ^{0,1} \mathcal{V}^*$. It follows from the results in \cite{N} (or \cite{BM} if $n=1$) that  $(M,g,\hat{J})$ is isomorphic to the twistor space $(\mathcal{Z},h_{t_1},J_2)$ of some compact quaternionic K\"ahler manifold with positive scalar curvature. Now the definitions of $\hat{J}$ and $J_2$ show that $(M,g,J)$ is isomorphic to  $(\mathcal{Z},h_{t_1},J_1)$. \hfill $\Box$

\vspace{3mm}
From the results of the previous sections we obtain

\begin{cor}\label{cor401}
A manifold, which satisfies the assumptions of Theorem~\ref{thm1}, is compact, simply connected, has positive and $\nabla^a$-parallel Ricci tensor, $\nabla^a T^a =0$ and $Hol(\nabla^a) \subset \rho (Sp(n)U(1))$. If furthermore the quaternionic K\"ahler base is not a symmetric space of rank greater than one, then $Hol(\nabla^a) = \rho (Sp(n)U(1))$.
\end{cor}

\vspace{3mm}
Theorem~\ref{thm1} has also the following consequences:

\begin{itemize}
\item In each dimension the manifolds satisfying its conditions are finitely many since the same is true for the compact quaternionic K\"ahler manifolds with positive scalar curvature \cite{L,LS}.
\item The only known examples of compact quaternionic K\"ahler manifolds with positive scalar curvature are the Wolf spaces \cite{W,B}, which are symmetric. Their twistor spaces $(\mathcal{Z},h_{t_1},J_1)$ are homogeneous. In fact, $(\mathcal{Z},h_{t_1})$ is a naturally reductive homogeneous space with canonical connection $\nabla^{a,t_1}$ iff the base manifold is symmetric (this follows from the proof of Proposition~\ref{prop6} for example).
\item It is known that in dimensions 4, 8 and 12 there are no compact quaternionic K\"ahler manifolds with positive scalar curvature other than the Wolf spaces \cite{FK,H,PS,HH}. Hence the only manifolds of dimension 6, 10 and 14, which satisfy the conditions of Theorem~\ref{thm1}, are the twistor spaces of $S^4 \cong \mathbb{H}P^1$, $\mathbb{C}P^2$; $\mathbb{H}P^2$, $Gr_2 (\mathbb{C}^4)$, $G_2/SO(4)$; $\mathbb{H}P^3$, $Gr_2 (\mathbb{C}^5)$, $\tilde{Gr}_4 (\mathbb{R}^7)$.
\end{itemize}

\section{The local case}\label{sec6}

The goal of this section is to prove the local version of Theorem~\ref{thm1}.

\begin{theo}\label{thm2}
Let $(M,g,J)$ be a Hermitian manifold such that $Hol(\nabla^a) \subset U(m) \times U(1)$, $T^a \in  \Lambda ^{2,0} \mathcal{H}^* \otimes \Lambda ^{0,1} \mathcal{V}^* \oplus \Lambda ^{0,2} \mathcal{H}^* \otimes \Lambda ^{1,0} \mathcal{V}^*$ and $T^a$ is non-degenerate. Then $m=2n$ and locally $(M,g,J)$ is isomorphic to the twistor space $(\mathcal{Z},h_{t_1},J_1)$ of some $4n$-dimensional quaternionic K\"ahler manifold $M'$ with positive scalar curvature. In particular, $\nabla^a T^a =0$ and $Hol(\nabla^a) \subset \rho (Sp(n)U(1))$, with equality if $M'$ is not locally symmetric of rank greater than one. 
\end{theo}

We begin the proof with the following straightforward

\begin{lem}\label{lem1}
Let $U,V \in \Gamma (\mathcal{V})$ and $X,Y \in \Gamma (\mathcal{H})$. Then
\begin{eqnarray}
& & \nabla _U V = {\nabla^a} _U V  \in \mathcal{V}, \nonumber \\
& & \nabla_U X = {\nabla^a} _U X - \frac{1}{2} T^a(U,X) \in \mathcal{H}, \nonumber \\
& & h\nabla_X U = \frac{1}{2} T^a(U,X), \qquad v\nabla _X U = {\nabla^a} _X U, \nonumber \\
& & h\nabla_X Y = {\nabla^a} _X Y, \qquad v\nabla_X Y = - \frac{1}{2} T^a(X,Y) \nonumber
\end{eqnarray}
($h$ and $v$ are the projections on $\mathcal{H}$ and $\mathcal{V}$ respectively). In particular, $\mathcal{V}$ is a totally geodesic distribution and therefore integrable.
\end{lem}

This lemma implies that each point of $M$ has a neighbourhood of the form $M' \times F$, where the fibres $\{ p' \} \times F$ are integral manifolds for $\mathcal{V}$. We restrict our considerations to this neighbourhood and denote it again by $M$. Let $\pi : M \longrightarrow M'$ be the projection.

\begin{lem}\label{lem2}
$(\mathfrak{L}_U g) (X,Y) = 0$ for $U \in \mathcal{V}$ and $X,Y \in \mathcal{H}$. Hence there exists a Riemannian metric $g'$ on $M'$ such that $\pi : (M,g) \longrightarrow (M',g')$ is a Riemannian submersion.
\end{lem}

\noindent {\it Proof:}
The first claim follows from Lemma~\ref{lem1}, the second one from Proposition~\ref{propa16}. \hfill $\Box$

\vspace{3mm}
From Proposition~\ref{prop5} we know that $m=2n$, $Hol(\nabla^a) \subset \rho (Sp(n)U(1))$ and $T^a = \lambda T_0 + \bar{\lambda} \overline{T}_0$, where $\lambda \not =0$ is a constant. Thus we have a $\nabla^a$-parallel and compatible with $g_{|\mathcal{H}}$ quaternionic structure $Q = span \{ I_{|\mathcal{H}}, J_{|\mathcal{H}}, K_{|\mathcal{H}} \}$ on $\mathcal{H}$, where $K$ and $I$ are defined by \pref{102}.

\begin{lem}\label{lem3}
Let $U \in \mathcal{V}$. Then $(h\mathfrak{L}_U (I_{|\mathcal{H}}), h\mathfrak{L}_U (J_{|\mathcal{H}}), h\mathfrak{L}_U (K_{|\mathcal{H}})) = (I_{|\mathcal{H}}, J_{|\mathcal{H}}, K_{|\mathcal{H}}).A(U)$, where
$$A(U) =
\begin{pmatrix}
0 & -2{\rm Re}(\lambda e^{2n+1}(U)) & \frac{1}{2} {\rm Im}(e^{2n+1}({\nabla^a} _U e_{2n+1})) \\
2{\rm Re}(\lambda e^{2n+1}(U)) & 0 & -2{\rm Im}(\lambda e^{2n+1}(U)) \\
-\frac{1}{2} {\rm Im}(e^{2n+1}({\nabla^a} _U e_{2n+1})) & 2{\rm Im}(\lambda e^{2n+1}(U)) & 0
\end{pmatrix}
.$$
Hence $Q = span \{ I_{|\mathcal{H}}, J_{|\mathcal{H}}, K_{|\mathcal{H}} \}$ projects on a quaternionic structure $Q'$ on $M'$.
\end{lem}

\noindent {\it Proof:}
\hfill The \hfill first \hfill claim \hfill follows \hfill from \hfill the \hfill definition \hfill of \hfill $K$ \hfill and \hfill $I$. \hfill Thus \\
$h\mathfrak{L}_U (I_{|\mathcal{H}}), h\mathfrak{L}_U (J_{|\mathcal{H}}), h\mathfrak{L}_U (K_{|\mathcal{H}}) \in span \{ I_{|\mathcal{H}}, J_{|\mathcal{H}}, K_{|\mathcal{H}} \}$ and therefore $Q$ is projectable  by Proposition~\ref{propa9}. \hfill $\Box$

\vspace{3mm}
Obviously $Q'$ and $g'$ are compatible, i.e., they define an almost quaternionic Hermitian structure on $M'$.

\begin{lem}\label{lem4}
$\nabla^a$ projects on the Levi-Civita connection $\nabla'$ of $(M',g')$.
Hence $Q'$ is $\nabla'$-parallel.
\end{lem}

\noindent {\it Proof:}
Since $\pi$ is a Riemannian submersion, $\nabla$ projects on $\nabla'$, i.e., $h\nabla_{X^h} Y^h = (\nabla'_X Y)^h$. Now Lemma~\ref{lem1} implies $h{\nabla^a} _{X^h} Y^h = (\nabla'_X Y)^h$, which means that $\nabla^a$ also projects on $\nabla'$ (see the Appendix). Since $Q$ is $\nabla^a$-parallel, the second claim follows from Lemma~\ref{lem3} and  Corollary~\ref{cora5}. \hfill $\Box$

\vspace{3mm}
Thus $(M',g')$ is quaternionic K\"ahler if $n>1$. Now we need to compute its curvature $R'$ to see that it is self-dual and Einstein if $n=1$ and that the scalar curvature $s'$ is positive for all $n$. This follows from

\begin{lem}\label{lem5}
$R' = \frac{|T^a|^2}{48n} R_{\mathbb{H}P^n} + R'_{hyper}$, where  $R'_{hyper}$ is the projection of $R_{hyper}$ (the hyper-K\"ahler part of $R^a$) and therefore has the symmetries of a hyper-K\"ahler curvature tensor.
\end{lem}

\noindent {\it Proof:}
According to O'Neill formulae \cite{O,B}
\begin{eqnarray}
& & R'(X,Y,Z,W) = R(X^h,Y^h,Z^h,W^h) \nonumber \\
& & \qquad + g(A_{Y^h} Z^h,A_{X^h} W^h) - g(A_{X^h} Z^h,A_{Y^h} W^h) -2g(A_{X^h} Y^h,A_{Z^h} W^h), \nonumber
\end{eqnarray}
where
$$A_{X^h} Y^h = \frac{1}{2} v[X^h,Y^h] = \frac{1}{2} v\nabla_{X^h} Y^h.$$
Hence, by Lemma~\ref{lem1}, $A_{X^h} Y^h = -\frac{1}{2} T^a(X^h,Y^h)$ and from \pref{204} we obtain
$$R'(X,Y,Z,W) = R^a(X^h,Y^h,Z^h,W^h) - g(T^a(X^h,Y^h),T^a(Z^h,W^h)).$$
This and Corollary~\ref{cor8} yield

\vspace{3mm}
$\qquad R'(X,Y,Z,W) = \frac{|T^a|^2}{48n} R_{\mathbb{H}P^n}^{g_{|\mathcal{H}}}(X^h,Y^h,Z^h,W^h) + R_{hyper}(X^h,Y^h,Z^h,W^h).$ \hfill   $\Box$

\vspace{3mm}
So, up to now we showed that for all $n$ $(M',g')$ is quaternionic K\"ahler with positive scalar curvature $s' = \frac{(n+2)|T^a|^2}{3}$. Let $\mathcal{Z}$ be its twistor space.

For $p \in M$ let $f(p)$ be the projection of $J_{|\mathcal{H}_p}$, i.e., $f(p) = {\pi_*}_{|\mathcal{H}_p}(J_{|\mathcal{H}_p}) \in Q'_{\pi (p)}$. Hence $f(p) \in \mathcal{Z}$ and in this way we obtain a map $f:M \longrightarrow \mathcal{Z}$.

Now it is straightforward (but somewhat long) to see that $f$ gives the desired isomorphism  of $(M,g,J)$ and $(\mathcal{Z},h_{t_1},J_1)$. The main points are to prove that
\begin{enumerate}
\item $f$ maps $g_{|\mathcal{V}}$ on ${h_{t_1}}_{|\mathcal{V}}$ and $J_{|\mathcal{V}}$ on ${J_1}_{|\mathcal{V}}$, and
\item $f_* X^h = X^h _Z$ ($X^h$ is the horizontal lift on $M$ and $X^h _Z$ is the horizontal lift on $\mathcal{Z}$).
\end{enumerate}
This completes the proof of Theorem~\ref{thm2}. \hfill $\Box$

\vspace{3mm}
Using Proposition~\ref{prop4}, Corollary~\ref{cor4} and Corollary~\ref{cor5} we obtain the obvious counterpart of Theorem~\ref{thm2} for nearly K\"ahler manifolds:

\begin{cor}\label{cor9}
Let $(M,g,J)$ be a strict nearly K\"ahler manifold with $Hol(\nabla^a) \subset U(m) \times U(1)$ and $T^a \in  \Lambda ^{2,0} \mathcal{H}^* \otimes \Lambda ^{1,0} \mathcal{V}^* \oplus \Lambda ^{0,2} \mathcal{H}^* \otimes \Lambda ^{0,1} \mathcal{V}^*$. Then $m=2n$ and locally $(M,g,J)$ is isomorphic to the twistor space $(\mathcal{Z},h_{t_1},J_2)$ of some $4n$-dimensional quaternionic K\"ahler manifold $M'$ with positive scalar curvature. In particular, $Hol(\nabla^a) \subset \rho_2 (Sp(n)U(1))$, with equality if $M'$ is not locally symmetric of rank greater than one.
\end{cor}

\noindent {\bf Remark:} The important point in Corollary~\ref{cor9} is that the result is local. For complete manifolds this follows from the results in \cite{N,BM} used in the proof of Theorem~\ref{thm1}.

\appendix}

\section{Projections of tensors and connections by submersions}

Let $\pi:M \longrightarrow M'$ be a (surjective) submersion. Denote by $\mathcal{V}$ the {\it vertical distribution} on $M$, i.e., $\mathcal{V} = Ker \, \pi_*$. Assume that the fibres of $\pi$ are connected and a distribution $\mathcal{H}$ is given on $M$, so that $TM=\mathcal{V} \oplus \mathcal{H}$. We call $\mathcal{H}$ {\it the horizontal distribution}.

Let  $v:TM \longrightarrow \mathcal{V}$ and $h:TM \longrightarrow \mathcal{H}$ be the projections. Since
$$T^r_s M = (TM)^{\otimes r} \otimes (T^*M)^{\otimes s}=(\mathcal{V} \oplus \mathcal{H})^{\otimes r} \otimes (\mathcal{V}^* \oplus \mathcal{H}^*)^{\otimes s} \supset \mathcal{H}^{\otimes r} \otimes  (\mathcal{H}^*)^{\otimes s}=T^r_s (\mathcal{H}),$$
we have a projection $h:T^r_s M \longrightarrow T^r_s (\mathcal{H})$.

Since $\pi$ is a submersion, $\pi_* \vert _{\mathcal{H}_p} : \mathcal{H}_p \longrightarrow T_{\pi (p)} M'$ is an isomorphism and defines an isomorphism of $T^r_s (\mathcal{H}_p)$ and $T^r_s (T_{\pi (p)} M')$.

\vspace{3mm}
\noindent {\bf Definition:} 1) Let $T \in T^r_s (T_p M)$. The tensor $\pi_* \vert _{\mathcal{H}_p} (hT)$ is called {\it the projection} of $T$.

2) Let $T' \in T^r_s (T_{\pi (p)} M')$. The tensor $T'^h:= (\pi_* \vert _{\mathcal{H}_p})^{-1}(T') \in T^r_s (\mathcal{H}_p) \subset T^r_s (T_p M)$ is called {\it the horizontal lift} of $T'$ in $p$.

3) Let $T' \in \Gamma (T^r_s M')$. Then by 2) $T'$ defines a tensor field $T'^h$, {\it the horizontal lift} of $T'$.

4) A tensor field $T \in \Gamma (T^r_s M)$ is said to be {\it projectable} if $hT=T'^h$ for some $T' \in \Gamma (T^r_s M')$. In this case $T'$ is called {\it the projection} of $T$.

\vspace{3mm}
The following proposition is straightforward.

\begin{prop}\label{propa1}
1) $T'^h (\omega_1 ^h,\dots,\omega_r ^h,X_1^h,\dots,X_s^h)= T' (\omega_1,\dots,\omega_r,X_1,\dots,X_s) \circ \pi$.

2) $\omega^h = \pi^* \omega$ for $\omega \in \Gamma (T^* M')$.

3) If $U \in \Gamma (\mathcal{V})$, then ${\mathfrak L}_U X^h \in \Gamma (\mathcal{V})$, ${\mathfrak L}_U (\pi^* \omega) =0$.

4) If $U \in \Gamma (\mathcal{V})$, $T \in \Gamma (T^r_s M)$,  then
$${\mathfrak L}_U (hT) (\omega_1 ^h,\dots,\omega_r ^h,X_1^h,\dots,X_s^h) = U(T (\omega_1 ^h,\dots,\omega_r ^h,X_1^h,\dots,X_s^h)).$$
In particular, $h{\mathfrak L}_U (hT)$ is linear with respect to $U$ and therfore depends pointwise on $U$.

5) $h{\mathfrak L}_U {T'}^h = 0.$
\end{prop}

\begin{prop}\label{propa6}
Let $T \in \Gamma (T^r_s M)$. Then $T$ is projectable iff $h{\mathfrak L}_U (hT) =0$ for each vertical $U$.
\end{prop}

\noindent {\it Proof:} $T$ is projectable iff $T (\omega_1 ^h,\dots,\omega_r ^h,X_1^h,\dots,X_s^h)$ is constant along the fibres for each $\omega_1,\dots,\omega_r,X_1,\dots,X_s$ on $M'$. Since we assume that the fibres are connected, this is equivalent to $U(T (\omega_1 ^h,\dots,\omega_r ^h,X_1^h,\dots,X_s^h)) =0$ for each vertical $U$. Thus the assertion follows from 4) in Proposition~\ref{propa1}. \hfill $\Box$

\vspace{3mm}
In a similar way we get

\begin{prop}\label{propa7}
Let $T \in \Gamma (T^r_s M)$ has no components in $\sum _{i=1}^s A_i$, where $A_i = \mathcal{H} \otimes \dots \otimes \mathcal{H} \otimes B_1^* \otimes \dots \otimes B_s^*$ with $B_i = \mathcal{V}$ and $B_j=\mathcal{H}$ for $j \not = i$. Then for each vertical $U$
$${\mathfrak L}_U T (\omega_1 ^h,\dots,\omega_r ^h,X_1^h,\dots,X_s^h) = U(T (\omega_1 ^h,\dots,\omega_r ^h,X_1^h,\dots,X_s^h)).$$
In particular, $T$ is projectable iff $h{\mathfrak L}_U T =0$ for each vertical $U$.
\end{prop}

Proposition~\ref{propa6} obviously follows from Proposition~\ref{propa7} since $T$ is projectable iff $hT$ is projectable.

\vspace{3mm}
\noindent {\bf Definition:} A subbundle $S$ of $T^r_s M$ is said to be {\it projectable} if there exists a subbundle $S'$ of $T^r_s M'$ such that
$\pi_* \vert _{\mathcal{H}_p} \circ h$ maps $S_p$ isomorphically onto $S'_{\pi (p)}$ for each $p \in M$ (in particular, $h:S_p \longrightarrow h(S_p)$ must be an isomorphism). $S'$ is called {\it the projection} of $S$.

\begin{prop}\label{propa8}
Let $S$ be a subbundle of $T^r_s M$ such that $h:S_p \longrightarrow h(S_p)$ is an isomorphism for each $p \in M$. Then $S$ is projectable iff $h{\mathfrak L}_U (hT) \in  \Gamma (hS)$ for each $T \in \Gamma (S)$ and each vertical $U$.
\end{prop}

\noindent {\it Proof:} The necessity is clear from 5) in Proposition~\ref{propa1}. For the sufficiency we have to prove that $\pi_* \vert _{\mathcal{H}_p} (h S_p) = \pi_* \vert _{\mathcal{H}_q} (h S_q)$ if $\pi (p) = \pi (q)$. We have the following straightforward

\begin{lem}\label{lema1}
Let $T_1, \dots , T_k$ be a local basis of $S$, $T'_1, \dots , T'_m$ be a local basis of $T^r_s M'$ and $U \in \mathcal{V}$. Let
$$(h{\mathfrak L}_U (hT_1), \dots , h{\mathfrak L}_U (hT_k)) = (hT_1, \dots , hT_k).A(U), \quad (hT_1, \dots , hT_k)=({T'_1}^h, \dots , {T'_m}^h).B,$$
where $A$ is a 1-form whose values are $k \times k$ matrices and $B$ is a function whose values are  $m \times k$ matrices. Then
$$(h{\mathfrak L}_U (hT_1), \dots , h{\mathfrak L}_U (hT_k)) = ({T'_1}^h, \dots , {T'_m}^h).dB(U)$$
and therfore
\begin{equation}\label{3}
dB(U)=B.A(U).
\end{equation}
\end{lem}

Since \pref{3} is a linear system for $B$, the span of the columns of $B(p)$ is independent of $p$. But this immediately implies
$$span \{ \pi_* \vert _{\mathcal{H}_p} (hT_1), \dots , \pi_* \vert _{\mathcal{H}_p} (hT_k) \} = span \{ \pi_* \vert _{\mathcal{H}_q} (hT_1), \dots , \pi_* \vert _{\mathcal{H}_p} (hT_k) \}$$
if $\pi (p) = \pi (q)$. \hfill $\Box$

\vspace{3mm}
In a similar way we have

\begin{prop}\label{propa9}
Let $S$ be a subbundle of $T^r_s M$ such that $h:S_p \longrightarrow h(S_p)$ is an isomorphism for each $p \in M$ and $S_p \cap \sum _{i=1}^s A_i = \{0\}$, where $A_i$ are the same as in Proposition~\ref{propa7}. Then $S$ is projectable iff $h{\mathfrak L}_U T \in  \Gamma (hS)$ for each $T \in \Gamma (S)$ and each vertical $U$.
\end{prop}

\noindent {\bf Definition:} 1) A connection $\nabla$ on $M$ is said to be {\it projectable} if there exists a connection $\nabla'$ on $M'$ such that $(\nabla'_X Y)^h = h\nabla_{X^h} Y^h$. $\nabla'$ is called {\it the projection} of $\nabla$.

2) For a connection $\nabla$ on $M$ we define $h\nabla$ by $(h\nabla)_X Y = h\nabla_{hX} (hY)$ and ${\mathfrak L}_X (h\nabla)$ by $({\mathfrak L}_X (h\nabla))_Y Z = {\mathfrak L}_X ((h\nabla)_Y Z)-(h\nabla)_{{\mathfrak L}_X Y} Z -(h\nabla)_Y ({\mathfrak L}_X Z)$.

\begin{prop}\label{propa110}
\hfill Let \hfill $\nabla$ \hfill be \hfill a \hfill connection \hfill on \hfill $M$ \hfill and \hfill $U$ \hfill be \hfill vertical. \hfill Then  \\
$h({\mathfrak L}_U (h\nabla))_{X^h} Y^h = h[U,h\nabla_{X^h} Y^h ]$.
\end{prop}

Using this, it is straightforward to prove

\begin{prop}\label{propa10}
A connection $\nabla$ on $M$ is projectable iff $h({\mathfrak L}_U (h\nabla))_{X^h} Y^h =0$ for each vertical $U$.
\end{prop}

\begin{prop}\label{propa11}
Let $\nabla$ project on $\nabla'$ and  $T \in \Gamma (T^r_s M)$  project on $T' \in \Gamma (T^r_s M')$. Then $\nabla (hT)$ projects on $\nabla' T'$.
\end{prop}

\begin{cor}\label{cora1}
Let $\nabla$ project on $\nabla'$ and  $T \in \Gamma (T^r_s M)$  project on $T' \in \Gamma (T^r_s M')$. Then $\nabla' T' =0$ iff $h\nabla (hT) =0$.
\end{cor}

\begin{prop}\label{propa12}
Let $\nabla$ project on $\nabla'$ and  $T \in \Gamma (T^r_s M)$  project on $T' \in \Gamma (T^r_s M')$. Let $T$ have no components in $A_1 \otimes \dots \otimes A_r \otimes A_{r+1}^* \otimes \dots \otimes A_{r+s}^*$, where one $A_i$ is $\mathcal{V}$ and the others are $\mathcal{H}$. Then $\nabla T$ projects on $\nabla' T'$.
\end{prop}

\noindent {\it Proof:} Under the given assumption $h\nabla T = h\nabla (hT)$ and the claim follows from Proposition~\ref{propa11}. \hfill $\Box$

\begin{cor}\label{cora2}
If $\nabla$, $\nabla'$, $T$, $T'$ are as in Proposition~\ref{propa12}, then $\nabla' T' =0$ iff $h\nabla T =0$.
\end{cor}

In a similar way as in Proposition~\ref{propa12} one proves

\begin{prop}\label{propa1000}
Let $\nabla$ project on $\nabla'$ and  $T \in \Gamma (T^r_s M)$  project on $T' \in \Gamma (T^r_s M')$. Let the decomposition $TM= \mathcal{V} \oplus \mathcal{H}$ be parallel with respect to $\nabla$. Then $\nabla T$ projects on $\nabla' T'$.
\end{prop}

\begin{prop}\label{propa13}
Let $\nabla$ project on $\nabla'$ and  $S$ be a subbundle of $T^r_s M$, which projects on $S'$. Then $S'$ is parallel with respect to $\nabla'$ iff $h\nabla _X (hT) \in h(S)$ for each $X \in \mathcal{H}$ and $T \in \Gamma (S)$.
\end{prop}

\begin{prop}\label{propa14}
Let $\nabla$ project on $\nabla'$ and  $S$ be a subbundle of $T^r_s M$, which projects on $S'$. Let $S \cap A_1 \otimes \dots \otimes A_r \otimes A_{r+1}^* \otimes \dots \otimes A_{r+s}^* = {0}$, where $A_i$ are as in  Proposition~\ref{propa12}. Then $S'$ is parallel with respect to $\nabla'$ iff $h\nabla _X T \in h(S)$ for each $X \in \mathcal{H}$ and $T \in \Gamma (S)$.
\end{prop}

\begin{cor}\label{cora5}
If $\nabla$, $\nabla'$, $S$, $S'$ are as in Proposition~\ref{propa14} and $\nabla_X S \subset S$ for each  horizontal X, then $S'$ is parallel with respect to $\nabla'$. This is true, in particular, if $S$ is parallel with respect to $\nabla$.
\end{cor}

\begin{prop}\label{propa15}
Let $\nabla$ project on $\nabla'$. Then the torsion $T$ of $\nabla$ projects on the torsion $T'$ of $\nabla'$.
\end{prop}

\begin{cor}\label{cora6}
If $\nabla$ projects on $\nabla'$ and $hT =0$, then $T' =0$. In particular, if $T =0$, then $T' =0$
\end{cor}

If $\pi:M \longrightarrow M'$ is a submersion and $g$ is a Riemannian metric on $M$, then we have a canonical horizontal distribution $\mathcal{H}=\mathcal{V}^\bot$.

\vspace{3mm}
\noindent {\bf Definition:} Let $\pi:(M,g) \longrightarrow (M',g')$ be a submersion. $\pi$ is called a {\it Riemannian submersion} if $g$ projects on $g'$.

\vspace{3mm}
As a consequence of Proposition~\ref{propa6} and Proposition~\ref{propa7}  we get

\begin{prop}\label{propa16}
Let $\pi:M \longrightarrow M'$ be a surjective submersion and $g$ be a Riemannian metric on $M$. Then the following conditions are equivalent:

1) There exists a Riemannian metric $g'$ on $M'$ such that $\pi:(M,g) \longrightarrow (M',g')$ is a Riemannian submersion.

2) $h{\mathfrak L}_U (hg) =0$ for each vertical U.

3) $h{\mathfrak L}_U g =0$ for each vertical U (this is true, in particular, if the vertical distribution has a basis of Killing vector fields).
\end{prop}

The following result is well known.

\begin{cor}\label{cora8}
If $\pi:(M,g) \longrightarrow (M',g')$ is a Riemannian submersion and $\nabla$ and $\nabla'$ are the Levi-Civita connections of $g$ and $g'$, then $\nabla$ projects on $\nabla'$.
\end{cor}

Let $\pi:M \longrightarrow M'$ be a submersion and $c$ be a conformal class on $M$. Then again we have a canonical horizontal distribution $\mathcal{H}=\mathcal{V}^\bot$.

\vspace{3mm}
\noindent {\bf Definition:} Let $\pi:(M,c) \longrightarrow (M',c')$ be a submersion. $\pi$ is called a {\it conformal submersion} if for each $p \in m$ and each $g \in c_p$ $g$ projects on an element $g' \in c'_{\pi (p)}$.

\vspace{3mm}
Obviously every conformal class $c$ on a manifold $M$ defines a 1-dimensional subbundle $S \subset T^0_2 M$.

\begin{prop}\label{propa17}
A submersion $\pi:(M,c) \longrightarrow (M',c')$ is a conformal submersion iff the subbundle $S$ defined by $c$ projects on the subbundle $S'$ defined by $c'$.
\end{prop}

\begin{prop}\label{propa18}
Let $\pi:M \longrightarrow M'$ be a surjective submersion and $c$ be a conformal class on $M$. Then the following conditions are equivalent:

1) There exists a conformal class $c'$ on $M'$ such that $\pi:(M,c) \longrightarrow (M',c')$ is a conformal submersion.

2) For each $g \in c$ there exists a 1-form $\alpha$ such that $h{\mathfrak L}_U (hg) =\alpha (U)hg$ for each vertical $U$.

3) There exist $g \in c$ and a 1-form $\alpha$ such that $h{\mathfrak L}_U (hg) =\alpha (U)hg$ for each vertical $U$.

4) For each $g \in c$ there exists a 1-form $\alpha$ such that $h{\mathfrak L}_U g =\alpha (U)g$ for each vertical $U$.

5) There exist $g \in c$ and a 1-form $\alpha$ such that $h{\mathfrak L}_U g =\alpha (U)g$ for each vertical $U$.
\end{prop}

From Corollary~\ref{cora6} we obtain

\begin{cor}\label{cora9}
If $\pi:(M,c) \longrightarrow (M',c')$ is a conformal submersion and $\nabla$ is a Weyl connection on $(M,c)$ which projects on a connection $\nabla'$, then $\nabla'$ is a Weyl connection on $(M',c')$.
\end{cor}

\vspace{3mm}
\noindent {\it Acknowledgement.} The author would like to thank Thomas Friedrich for
drawing his attention to the problem considered in the present paper and for many useful discussions.

\vspace{10mm}
\noindent
Bogdan Alexandrov \\
Universit\"at Greifswald \\
Institut f\"ur Mathemathik und Informatik \\
Friedrich-Ludwig-Jahn-Stra{\ss}e 15a \\
17487 Greifswald \\
{\tt e-mail: \quad boalexan@uni-greifswald.de}

\end{document}